\documentclass{article}

\newcommand{\R}{I\!\! R}

\newcommand{\Z}{I\!\! Z}

\begin{document}

\centerline  {\bf Fibr\'es
affines.}

\medskip

par Tsemo Aristide,

The Abdus Salam Centre for Theoretical Physics

Strada Costiera, 11.

Trieste, Italy.

\bigskip

\centerline{\bf
Abstract.}

{\it We study affine maps between affine manifolds. Even
when two fibers are compact and diffeomorphic, they can
inherit different affine structures from the source space. This leads
to a fixed linear holonomy deformation theory of the affine structure
of an affine manifold. We found various conditions which make the
fibers to be affinely isomorphic. We  also classify affine bundles
which total space is a compact and complete affine
manifold.}

\bigskip

\centerline{\bf R\'esum\'e.}

{\it On
\'etudie les applications affines entre vari\'et\'es affines. Deux
fibres distinctes, compactes, et diff\'eomorphes peuvent h\'eriter de
structures affines distinctes de la vari\'et\'e source. Ceci conduit
\`a une th\'eorie de d\'eformation  \`a holonomie lin\'eaire fix\'ee
de la structure affine d'une vari\'et\'e affine. On trouve
plusieurs conditions sous lesquelles les fibres sont isomorphes deux
\`a deux.  On classifie les fibr\'es affines  dont l'espace total est
une vari\'et\'e affine compacte et compl\`ete.}

\bigskip

{\bf
Introduction.}

\bigskip

On \'etudie  les applications affines
entre vari\'et\'es affines, en particulier  les fibr\'es affines: ce
sont les applications affines surjectives.  La classification des
applications affines est envisageable si on suppose que  les fibres
sont au moins hom\'eomorphes entre elles. C'est par exemple le cas  si
l'espace source de  l'application est une vari\'et\'e affine compacte.
 Voici un exemple qui montre que la situation peut \^etre beaucoup
plus compliqu\'ee. Consid\'erons une surjection de ${\R}^m$ dans
${\R}^p$, sur chaque fibre  de la surjection, enlevons un
sous-ensemble, on obtient ainsi une application affine, telle que les
structures topologiques de deux  fibres distinctes  peuvent
\^etre tr\`es diff\'erentes. Il est donc difficile d'aborder le
probl\`eme de classification dans  toute sa g\'en\'eralit\'e, car on ne
peut pas pr\'eciser les donn\'ees \`a partir desquelles se fera la
classification. Dans les cas que nous consid\'erons, les fibres d'un
fibr\'e affine sont des vari\'et\'es diff\'erentiables, deux \`a deux
hom\'eomorphes, dont les structures affines peuvent \^etre
diff\'erentes bien qu'elles aient la m\^eme holonomie lin\'eaire.

Si le second groupe d'homotopie de la base d'un fibr\'e affine \`a
espace total compact est nul (c'est le cas si elle est
compl\`ete), la suite exacte de Serre montre que le groupe
fondamental de l'espace total $\pi_1(M)$, est une extension de
celui de la base $\pi_1(B)$, par celui de la fibre $\pi_1(F)$.

Posons ${\R}^{m+l}={\R}^m\oplus
{\R}^l$, et notons $Aff({\R}^{m},{\R}^l)$ le sous-groupe
des automorphismes affines de ${\R}^{m+l}$ dont la partie
lin\'eaire fixe ${\R}^l$,
 et $s$ la surjection canonique de
${\R}^{m+l}$ dans ${\R}^m$.
 Il existe une application:
 $s':Aff({\R}^m,{\R}^l)\rightarrow Aff({\R}^m)$
 $\gamma\mapsto s'(\gamma)(x):=s(\gamma(x,0))$.
 Notons $Aff_I({\R}^m,{\R}^l)$ le sous-groupe de $Aff({\R}^m,{\R}^l)$
 qui se projette sur ${\R}^m$ en l'identit\'e.

  On en d\'eduit le probl\`eme
alg\'ebrique suivant analogue au notre: \'etant donn\'es deux groupes
$\pi_1(B)$ et $\pi_1(F)$, deux repr\'esentations
affines respectivement de $\pi_1(B)$ dans $Aff({\R}^m)$, et
de $\pi_1(F)$ dans $Aff_I({\R}^m,{\R}^l)$.   Classifier tous les
diagrammes commutatifs suivants $$ \matrix {& 1  \rightarrow &
\pi_1(F) &\rightarrow & \pi_1(M)  \rightarrow&\pi_1(B)&
\rightarrow  1 \cr& \ \  & \downarrow &\ \ & \downarrow &\ \
\downarrow \cr & 1\rightarrow &Aff_I({\R}^m,{\R}^l) &\rightarrow
&Aff({\R}^{m},{\R}^l)& \rightarrow Aff({\R}^m)& \rightarrow 1},$$

o\`u la premi\`ere  suite horizontale est
exacte, sous la relation d'\'equivalence suivante: Soient
$1\rightarrow\pi_1(F)\rightarrow\pi_1(M_{i})\rightarrow
\pi_1(B)\rightarrow 1$, $i=\{1,2\}$ deux suites exactes associ\'ees
\`a deux diagrammes du type pr\'ec\'edent. Ces diagrammes sont
\'equivalents si et seulement si les extensions sont isomorphes et il
existe une transformation affine de $Aff({\R}^{m},{\R}^l)$ qui
conjugue la repr\'esentation $\pi_1(M_1)\rightarrow
Aff({\R}^{m},{\R}^l)$ en la
repr\'esentation $\pi_1(M_2)\rightarrow Aff({\R}^{m},{\R}^l)$.

 La donn\'ee d'un fibr\'e affine induit une
d\'eformation de la structure affine d'une fibre param\'etr\'ee par la
base. Soient $\Gamma$ un groupe de type fini, et $E$ un espace
vectoriel muni d'une action lin\'eaire de $\Gamma$.  L'approche de ce
probl\`eme pose celui de d\'eterminer la structure de l'ensemble
(ou d'un de ses ouverts (voir [G3] p.178)) des
orbites des \'el\'ements de $H^1(\Gamma,E)$ sous l'action d'un groupe
de gauge. Des hypoth\`eses faites sur la structure de cet espace des
modules par exemple, qu'il  est une vari\'et\'e alg\'ebrique
rendant alg\'ebrique la projection de $H^1(\Gamma,E)$ sur lui,
entraine que les fibres sont des vari\'et\'es affines isomorphes entre
elles si la base est une vari\'et\'e affine compacte et compl\`ete.

- On
d\'emontre aussi que les fibres d'un fibr\'e affine dont l'espace total
est une vari\'et\'e affine compacte et compl\`ete, et le groupe
fondamental d'une fibre est nilpotent sont isomorphes entre elles.

\medskip

-{\bf On classifie} les fibr\'es affines \`a espace total compact dont
la d\'eveloppante est injective et tels que le fibr\'e relev\'e sur le
rev\^etement universel de la base est un produit affine,

- Plus g\'en\'eralement, on classifie les fibr\'es affines \`a espace
total compact et complet.

Les invariants d\'etermin\'es dans le
premier cas sont plus intrins\`eques \`a la g\'eom\'etrie affine que
ceux du second cas.

\medskip

Dans [T5] on donne une nouvelle classification des fibr\'es affines
\`a espace total compact et complet en utilisant la th\'eorie des gerbes.

 Une question importante est de savoir
dans quelle mesure la classe d'isomorphisme d'un fibr\'e affine peut
d\'eterminer celle de son espace total en tant que vari\'et\'e
affine. Soient $f_1$ et $f_2$ deux fibr\'es affines  d'espace total
compact (de m\^eme dimension) au-dessus d'une m\^eme
base $(B,\nabla_B)$. On suppose que $(B,\nabla_B)$  ne  peut avoir de
feuilletage affine \`a feuilles compactes, et en outre que la dimension
des fibres  de $f_1$ et $f_2$ est strictement inf\'erieure \`a celle
de $B$. Il est facile de prouver que tout isomorphisme affine  entre
les espaces total de $f_1$ et $f_2$, est un isomorphisme de fibr\'es
affines. Les  classes d'isomorphismes de ces fibr\'es affines,
d\'eterminent les classes d'isomorphismes de leurs espaces total en
tant que vari\'et\'es affines. Les vari\'et\'es affines de
dimension $3$ dont le groupe fondamental n'est pas polycyclique, et
les vari\'et\'es de Hopf   ne peuvent avoir de feuilletages affines
\`a feuilles compactes.

 \medskip

  Voici le plan de notre
travail:

0. Introduction

 1. G\'en\'eralit\'es

 2. Adh\'erence de
Zariski et structure des fibres.

3. D\'eformation des structures
affines \`a holonomie lin\'eaire fix\'ee.

4. Les fibr\'es affines
affinement localement triviaux.

5. Le cas g\'en\'eral.

\medskip

{\bf 1. G\'en\'eralit\'es.}

\medskip

Une vari\'et\'e affine
$(M,\nabla_M)$ de dimension $n$, est une
vari\'et\'e diff\'erentiable $M$, de dimension $n$, munie d'une
connexion  $\nabla_M$ dont les tenseurs de courbure et de torsion sont
nuls. La connexion $\nabla_M$ d\'efinit sur $M$ un atlas (affine) dont
les fonctions de transitions sont des restrictions d'\'el\'ements de
$Aff({\R}^n)$. Soient $(M,\nabla_M)$ et $(M',\nabla_{M'})$
deux vari\'et\'es affines dont les structures affines sont
respectivement d\'efinies par les atlas affines $(U_i,\phi_i)$ et
$(U'_j,\phi'_{j'}).$

 Une application  diff\'erentiable $f$, de
$(M,\nabla_M)$ dans $(M',\nabla_{M'})$, est affine si et seulement
si ${\phi'_{j'}} \circ f_{\mid
U_i} \circ{\phi_i}^{-1}$ est affine.  On
note $App((M,\nabla_M),(M',\nabla'_{M'}))$ l'ensemble des applications
affines de $(M,\nabla_M)$ dans $(M',\nabla'_{M'})$, et
$Aff(M,\nabla_M)$ l'ensemble des automorphismes affines
de $(M,\nabla_M)$.

 La structure affine de  $(M,\nabla_M)$ se
rel\`eve sur son rev\^etement universel $\hat M$,  en une structure
affine $(\hat M, \hat \nabla_{\hat M})$ pour laquelle la projection
rev\^etement  $p_M: (\hat M,\hat \nabla_{\hat M}) \rightarrow
(M,\nabla_M)$ est une application affine. La structure affine
de $\hat M$ est d\'efinie par un diff\'eomorphisme local  affine
$D_M: (\hat M,\hat\nabla_{\hat M}) \rightarrow ({\R}^n,\nabla_0)$,
o\`u $\nabla_0$ est la connexion standard de ${\R}^n$ d\'efinie
par $$\nabla_XY = DY(X)$$  $X$,  $Y$ sont deux champs de vecteurs de
${\R}^n$, et $DY$  la diff\'erentielle de $Y$.

 La d\'eveloppante
donne lieu \`a une repr\'esentation $A_M: Aff(\hat M,\hat\nabla_{\hat
M})\rightarrow Aff({\R}^n)$ d\'efinie par le diagramme
commutatif suivant: $$ \matrix
{(\hat M,\hat\nabla_{\hat M})&{\buildrel{g}\over{\longrightarrow}}&
(\hat M,\hat\nabla_{\hat M})\cr\downarrow D_M &\ \ \ \ & \downarrow
D_M\cr {\R}^n &{\buildrel{A_M(g)}\over{\longrightarrow}}&
{\R}^n}$$  o\`u $g$ est un \'el\'ement de
$Aff(\hat M,\hat\nabla_{\hat M})$. La restriction de $A_M$, au groupe
fondamental $\pi_1(M)$, de $M$, est la
repr\'esentation d'holonomie $h_M$. La partie lin\'eaire $L(h_M)$, de
$h_M$, est la repr\'esentation d'holonomie lin\'eaire de
 $(M,\nabla_M)$.

 \medskip

 Soit $f$ un \'el\'ement de
$App((M,\nabla_M),(M',\nabla'_{M'}))$, $f$ se rel\`eve en une
application affine  $\hat f: (\hat M,\hat \nabla_{\hat M})\rightarrow
(\hat M',\hat \nabla'_{\hat M'})$. Soit $U$ un ouvert connexe de
$\hat M$, sur lequel la restriction de la d\'eveloppante est
injective, $D_{M'{_{\mid \hat f(U)}}}\circ \hat f_{\mid U}\circ
{D_M}^{-1}_{\mid D_M(U)}$ est une application affine de
$D_M(U)\rightarrow {\R}^{n'}$, qui se prolonge en une
application affine  $f_1:{\R}^n\rightarrow {\R}^{n'}$. Celle-ci v\'erifie
$$D_{M'}\circ \hat f = f_1\circ
D_M.$$ L'application $\hat f$ donne lieu \`a un morphisme de
groupes $\hat f_{\pi}:\pi_1(M)\rightarrow \pi_1(M')$ d\'efini par la
relation d'\'equivariance $$\hat f\circ\gamma =\hat
f_{\pi}(\gamma)\circ \hat f.$$  On en d\'eduit un morphisme
$f_{\pi}: h_M(\pi_1(M))\rightarrow h_{M'}(\pi_1(M'))$ tel que $$
f_1\circ \gamma = f_{\pi}(\gamma)\circ f_1\leqno (0.1)$$ si $h_M$ est
injective. L'application $f_1$ transforme $D_M(\hat M)$ en un
sous-ensemble de $D_{M'}(\hat M')$.

\medskip

{\bf Question}:
Etant donn\'ees $f_{\pi}$ et $f_1$ telles que $f_1$ envoie $D_M(\hat
M)$ dans $D_{M'}(\hat M')$, et $f_1\circ h_M(\gamma) =
f_{\pi}(h_M(\gamma))\circ f_1$. A quelles conditions
peut-on reconstruire $f$?

 La r\'eponse \`a la question
pr\'ec\'edente est positive si  $D_{M'}$ est
injective.

 L'exemple suivant montre que si on se donne $f_1$ et
$f_{\pi}$ v\'erifiant $(0.1)$, il n'existe pas toujours d'application
$f$.

 Munissons $T^3$ de la structure affine de Sullivan-Thurston [S-T
p.22]. L'image de la d\'eveloppante est ${\R}^3=Vect(e_1,e_2,e_3)$
 priv\'e des axes de coordonn\'ees,
et l'holonomie est engendr\'ee par deux applications diagonales dans
ces axes de coordonn\'ees ayant des valeurs propres strictement
positives. Le quotient de $Vect(e_1,e_2)- \{{\R}e_1\cup {\R}e_2\}$ par
 la restriction de l'holonomie est la r\'eunion de $4$
tores.  On peut  choisir pour $f_1$ la projection de ${\R}^3$ sur
$Vect(e_1,e_2)$ parall\`element \`a $Vect\{e_3\}$, et $f_{\pi}$ la
restriction de l'holonomie \`a $Vect(e_1,e_2)$. Il n'existe pas
d'application $f$ correspondante. Dans ce cas $f_1$ n'est
pas d\'efinie sur l'image de la d\'eveloppante de la structure affine
de $T^3$.

\medskip

{\bf D\'efinition.}

- On appelle fibration
affine une application affine surjective.

- Un isomorphisme affine
entre deux fibr\'es affines de m\^eme base est un isomorphisme affine
entre leur espace total qui se projette sur leur base en une
transformation affine, autrement dit un isomorphisme qui envoie une
fibre sur une fibre.

- Un feuilletage affine sur une vari\'et\'e
affine $(M,\nabla_M)$ de dimension $n$, est un feuilletage de $M$
dont le relev\'e sur $\hat M$ est l'image r\'eciproque par $D_M$ de la
restriction \`a $D_M(\hat M)$ d'un feuilletage par
sous-espaces affines parall\`eles de ${\R}^n$.

\medskip

{\bf
Remarque.} Le produit fibr\'e de deux fibr\'es affines est un fibr\'e
affine.

\medskip

{\bf Proposition 1.1 [T4].} {\it Soit $f$ une
application affine, de source la vari\'et\'e affine compacte
et connexe $(M,\nabla_M)$, alors la distribution $F^f_x:=\{v\in T_xM
/df_x(v) = 0\}$ d\'efinit sur $M$ un feuilletage ${\cal F}^f$, dont les
feuilles sont les fibres d'une fibration affine.}

\medskip

La
base de la fibration affine de la proposition pr\'ec\'edente n'est pas
forc\'ement l'image de l'application. Si la source de $f$ est
compacte, l'image de $f$ est une sous-vari\'et\'e affine avec
\'eventuellement des points de self-intersections.

\medskip

{\bf Proposition 1.2.}
{\it  Supposons que $M$ est
compacte, alors l'image de $f$ a des points de self-intersections si
et seulement si la fonction qui a un point $x$ de $f(M)$, associe le
cardinal des composantes connexes de l'image r\'eciproque  de $x$ par
$f$  n'est pas constante. }

\medskip

{\bf Preuve.} Supposons que
$f(M)$ a un point de self-intersection $y$. Le point $y$ est
contenu dans un ouvert de carte affine $V'$, tel que l'ensemble des
points de self-intersections de $V'$ est une sous-vari\'et\'e
affine. Soit $x_i$ un \'el\'ement d'une composante connexe $C_i$,
de $f^{-1}(y)$. Il existe un ouvert connexe $U_i$ de $M$ contenant
$C_i$, tel que $f(U_i)$ est inclu dans $V'$, et  les ensembles $U_i$
sont deux \`a deux disjoints. Choisissons un ouvert $V$ inclu
dans $V'$ tel que $f^{-1}(V)$ est inclu dans l'union des $U_i$. Ceci
est possible car $M$ est compacte. Puisque $y$ est un point
de self-intersection, il existe au moins $i,j$ tels que
$dim(f(U_i)\cap f(U_j))<dimf(U_i)$. Les
ensembles $f^{-1}((f(U_i)-f(U_j))\cap V)$
et $f^{-1}((f(U_j)-f(U_i))\cap V)$ sont disjoints. On en d\'eduit que
pour tout point $y_i$ de $(f(U_i)-f(U_j))\cap V$, le cardinal
des composantes connexes de $f^{-1}(y_i)$, est strictement inf\'erieur
\`a celui  de celles de $f^{-1}(y)$.

 Pour tout point $y$ de
$f(M)$, on va noter $n_y$, le cardinal des composantes connexes de
l'image r\'eciproque de $y$ par $f$. Supposons maintenant que la
fonction $y\mapsto n_y$ n'est pas constante, et $f(M)$ n'a pas de
point de self-intersection. Soit $Z$ l'ensemble des points $z$ de
$f(M)$, tels que  $n_z$ est maximal. Si $Z$ ne contient pas un
point de self-intersection, alors $Z$ est ouvert dans $f(M)$. Or
$Z$ est ferm\'e car $M$ est compact. On en d\'eduit que $Z=f(M)$.
Donc la fonction $z\mapsto n_z$ est constante. Il y a
contradiction$\bullet$

\medskip

 Un classique th\'eor\`eme d'Ehresmann
[God], affirme  que toute submersion de source une vari\'et\'e
diff\'erentiable compacte est une fibration
diff\'erentiable localement triviale. Une fibration affine d'espace
total  compact $M$,  est donc une fibration localement triviale
pour la structure diff\'erentiable. On note $F$ sa fibre type.  Dans
la suite, on supposera que $f$ est localement triviale pour la
structure diff\'erentiable.  Le feuilletage ${\cal F}^f$  se
rel\`eve sur $\hat M $, en un feuilletage $\hat{\cal F}^f$. Ce
dernier est le relev\'e par $D_M$, d'un feuilletage $D_M(\hat {\cal
F}^f)$ de ${\R}^n$ par des sous-espaces
affines parall\`eles.

\medskip

Ecrivons ${\R}^n ={\R}^m\oplus {\R}^l$, o\`u ${\R}^m$ est un
 suppl\'ementaire de
la direction ${\R}^l$, des feuilles de $D_M(\hat {\cal
F}^f).$ Pour tout \'el\'ement $\gamma$ de $\pi_1(M)$, $$h_M(\gamma)
= (A_{\gamma}(x)+a_{\gamma},B_{\gamma}(y)+ C_{\gamma}(x)+d_{\gamma}).$$ Comme
le feuilletage ${\cal F}^f$ n'a pas d'holonomie,  pour tout \'el\'ement
$\gamma$ de $\pi_1(F)$, on a
$$h_M(\gamma)= (x,B_{\gamma}(y)+C_{\gamma}(x)+d_{\gamma}).$$ On en
d\'eduit que l'holonomie lin\'eaire d'une fibre ne d\'epend pas de
celle-ci.

\medskip

Pour tout \'el\'ement $x$ de ${\R}^m$, l'application

$$\pi_1(F)\longrightarrow {\R}^l$$
 $$\gamma\longmapsto C_\gamma(x)+d_\gamma$$

  est un $1-$cocycle
pour l'action de l'holonomie lin\'eaire d'une fibre. Notons $r(x)$ sa
classe de cohomologie. On en d\'eduit l'existence d'une application
affine $r:{\R}^m\rightarrow H^1(\pi_1(F),{\R}^l)$ qui a $x$
associe $r(x)$, $H^*(\pi_1(F),{\R}^l)$ est le $*$ groupe
de cohomologie de $\pi_1(F)$ relatif \`a l'action de l'holonomie
lin\'eaire d'une fibre. Si la base est compl\`ete, alors l'image de $r$
est un sous-espace affine.

Si les fibres sont compactes et
compl\`etes, alors l'image de $r$ est contenue dans
la sous-vari\'et\'e alg\'ebrique $L$, de $H^1(\pi_1(F),{\R}^l)$,
d\'efinie par $L=\{c\in H^1(\pi_1(F),{\R}^l)/\Lambda^l c\neq 0\}$
(Voir [F-G-H] th\'eor\`eme 2.2).

Il existe une vari\'et\'e
affine $(M,\nabla_M)$ compacte, incompl\`ete, et de dimension $3$,
telle que $\Lambda^3c_M\neq 0$ (Voir [G1]).

 L'exemple suivant
largement inspir\'e de [F-G-H] p. $510-511$, montre que deux
fibres distinctes d'un fibr\'e affine dont les fibres sont compactes
peuvent \^etre des vari\'et\'es affines qui ne sont pas
isomorphes. Soit $s$ un r\'eel, consid\'erons le sous-groupe
$G_s'$, de  $Aff({\R}^3)$, dont tout \'el\'ement $g$ est de la
forme $$g_{u,v,t,s}(x,y,z) =(e^{2t}x+uz,y+st+vz,e^{t}z),$$ o\`u
 $t,u,v\in {\R}$. Le groupe $G_s'$ pr\'eserve le demi-espace $H$ de
${\R}^3$  d\'efini par $z>0$. Il contient
un sous-groupe cocompact $G_s$. Le quotient de $H$ par $G_s$ est une
vari\'et\'e affine compacte $M_s$. Si $s$ d\'esigne un r\'eel non nul,
$M_s$ n'est pas radiante, alors que $M_0$ l'est. Donc les vari\'et\'es
affines $M_s$ et $M_0$ ne sont pas isomorphes (voir [F-G-H]).
On peut supposer que $G_s$ est l'image de $G_0$ par l'application
  $g_s:G_0\rightarrow G_s$, $g_{u,v,t,0}\mapsto
g_{u,v,t,s}$. Consid\'erons
maintenant la vari\'et\'e affine $M$, de dimension $4$, quotient de
 ${\R}\times H$ par l'action  de $G_0$ d\'efinie par
$g_{u,v,t,0}(s,y)=(s,g_s(g_{u,v,t,0})(y))$,  la projection de
 ${\R}\times H\rightarrow {\R}$, passe au quotient en une projection
de $M$ sur ${\R}$, on obtient ainsi un fibr\'e affine au-dessus
de ${\R}$, dont les fibres ne sont pas isomorphes entre
elles.

\medskip

Soit ${\R}_+^*$ l'ensemble des r\'eels strictement positifs.
Construisons maintenant un fibr\'e
affine au-dessus de ${\R}_+^*$, dont les fibres sont des tores
munis de structures affines deux \`a deux
distinctes.

\medskip

Soient $s$ un \'el\'ement de ${\R}^*_+$,
et $\phi_s$ l'application polynomiale de ${\R}^2$, d\'efinie
par
 $\phi_s(x,y)=(x+{1\over s}y^2,y)$.
  Conjuguons le groupe $\Gamma_s$
engendr\'e par
 $l_{1,s}(x,y)=(x+1,y)$
 et $l_{2,s}(x,y)=(x,y+s)$ par
$\phi_s$.

On obtient le groupe $\Gamma'_s$ engendr\'e
par
 $$\phi_{s}\l_{1,s}\phi_{s}^{-1}=l_{1,s},$$
  et
 $$l'_{2,s}= \phi_{s}l_{2,s}(\phi_{s})^{-1}(x,y)=(x+2y+s,y+s).$$
  Soient
$s$ et $s'$ deux r\'eels positifs non nuls, et $g_{s,s'}$ l'application
de ${\R}^2$ d\'efinie par
$g_{s,s'}(x,y)=(x,{s'\over s}y)$.

En
conjuguant $\Gamma'_s$ par $f_{s,s'}= \phi_{s'} g_{s,s'}{\phi_{s}}^{-1}$,
on obtient $\Gamma'_{s'}$. Remarquons que $f_{s,s'}$ n'est pas affine
si $s$ est diff\'erent de $s'$. Si les tores affines $T_s$, et
$T_{s'}$, quotients de ${\R}^2$ par $\Gamma'_s$, et
$\Gamma'_{s'}$,  \'etaient affinement isomorphes, il existerait
une transformation affine $A_{s',s}$ qui conjuguerait $\Gamma'_{s'}$
en $\Gamma'_{s}$. La partie lin\'eaire $L(A_{s',s})$ de cette
transformation serait de la forme
$$ L(A_{s',s})=\pmatrix{a & b\cr 0 &
a}$$

car la partie lin\'eaire du conjugu\'e de  ${l'}_{2,s'}$ par
$A_{s',s}$ est celle de ${l'}_{2,s}$.
En effet sans restreindre la g\'en\'eralit\'e, on peut supposer que
${l'}_{2,s'}$ est conjugu\'e en ${l'}_{2,s}$, puisqu'il existe un
 automorphisme
affine de $T_s$ qui transforme ${l'}_{2,s}$ en ${{l'}_{2,s}}^{-1}$.

On en d\'eduit que $a=1$, car
$l_{1,s}$ est conjugu\'e en $l_{1,s}$
par $A_{s',s}$. Il en r\'esulte
que la seconde coordonn\'ee du conjugu\'e de $l_{2,s'}$ par $A_{s',s}$
est $s'$ car elle ne peut \^etre n\'egative.
  Ceci n'est pas possible si $s$ est diff\'erent de $s'$. Les
quotients de ${\R}^2$ par $\Gamma'_s$ et $\Gamma'_{s'}$ sont
des vari\'et\'es affines qui ne sont pas isomorphes si $s$ est
distinct de $s'$.

 Consid\'erons maintenant les transformations affines
de ${\R}^*_+\times{\R}^2$ d\'efinies
par
 $\gamma_1(x,y,z)=(x,y+1,z)$
  et
  $\gamma_2(x,y,z)=(x,y+2z+x,z+x).$

  Le
quotient de ${\R}^*_+\times {\R}^2$ par le groupe engendr\'e
par $\gamma_1$ et $\gamma_2$,  est l'espace total d'un fibr\'e affine
au-dessus de ${\R}^*_+$, dont les fibres ont des structures affines
deux \`a deux distinctes.

 Plus tard (Proposition 4.1) on verra que
si la base d'un fibr\'e affine dont les fibres sont munies d'une
structure affine compl\`ete du tore de dimension $2$  distincte de
la structure riemannienne plate est compl\`ete, alors toutes les
fibres sont isomorphes entre elles.

\medskip

{\bf Question:} Les
fibres d'une fibration affine sont-elles isomorphes entre elles
si l'espace total est compact ?

 \medskip

Supposons que les fibres
sont compl\`etes,  et soient $x_0$ et $x_1$ deux points de la base
$(B,\nabla_B)$. Les fibres au-dessus de $x_0$, et $x_1$, sont
isomorphes, si et seulement si leurs holonomies respectives  sont
conjugu\'ees par une application affine. Les fibres d'un fibr\'e
affine sont deux \`a deux affinement isomorphes, si le groupe des
automorphismes affines qui pr\'eserve la fibration se projette sur la
base en un groupe qui agit transitivement su
celle-ci.

\medskip

{\bf Remarque.} Dans [G1], il est donn\'e
des exemples de structures de vari\'et\'es affines sur une vari\'et\'e
diff\'erentiable qui sont topologiquement \'equivalentes sans \^etre
diff\'erentiablement \'equivalentes.

\medskip

 De la suite
exacte d'homotopie associ\'ee \`a une fibration, on d\'eduit la suite
courte suivante: $ \pi_2(B)\rightarrow \pi_1(F) \rightarrow \pi_1(M)
\rightarrow \pi_1(B) \rightarrow 1$. Si  $\pi_2(B)=1$, (c'est le cas
si la base est compl\`ete) $\pi_1(M)$ est une extension de $\pi_1(B)$
par $\pi_1(F)$.

\medskip

{\bf Proposition 1.3.}
{\it  Supposons
que les fibres de $f$ sont compactes, et l'espace total du relev\'e
$\hat f$, de $f$ sur $\hat B$, a un feuilletage transverse aux
fibres, alors l'application $g:\pi_1(F)\rightarrow\pi_1(M)$, de la
suite exacte d'homotopie est injective.}

\medskip

{\bf
Preuve.} S'il existe un feuilletage transverse, aux fibres de $\hat
f$, alors d'apr\`es un th\'eor\`eme d'Ehresmann [God], $\hat f$ est
une suspension diff\'erentiable car ses fibres  sont compactes. Il en
r\'esulte que $\hat f$ est un produit diff\'erentiable, car $\hat B$
est simplement connexe, et par suite $\pi_1(F)$ s'injecte dans
$\pi_1(M)\bullet$

\medskip

On  d\'eduit que $\pi_1(M)$, est une
extension de $\pi_1(B)$ par $\pi_1(F)$, sous les conditions de la
proposition pr\'ec\'edente.

 Si la d\'eveloppante de l'espace total
est injective, $\pi_1(M)$ est une extension de $\pi_1(B)$ par
$\pi_1(F)$ car dans ce cas, $\pi_1(F)$ est le noyau de la restriction
de la repr\'esentation d'holonomie de l'espace total \`a ${\R}^m$.

\medskip

{\bf 2. Adh\'erence de Zariski et structures
des fibres.}

\medskip

 Dans ce paragraphe, on va
d\'eterminer diff\'erentes conditions, sous lesquelles les fibres d'un
fibr\'e affine sont deux \`a deux isomorphes. La majorit\'e de ces
conditions s'expriment en fonction d'adh\'erences de
Zariski.

\medskip

{\bf Proposition 2.1.} {\it Soit $f$ un
fibr\'e affine de base une vari\'et\'e affine compacte et
compl\`ete $(B,\nabla_B)$. Si $f$ a une fibre radiante, toutes ses
fibres sont radiantes. De plus si les fibres sont compactes, alors
leurs structures affines sont isomorphes.}

\medskip

{\bf
Preuve.} S'il existe une fibre radiante, alors il existe
un point $x_0$ de ${\R}^m$ tel que $r(x_0)=0$. On
a $r(\gamma(x_0))=0$, pour tout \'el\'ement $\gamma$ de $\pi_1(B)$.
L'ensemble $\{\gamma(x_0)$, $\gamma\in \pi_1(B)\}$ engendre l'espace
affine ${\R}^m$ (voir [F-G-H] th\'eor\`eme 2.2). On en d\'eduit
que $r=0$, car $r$ est une application affine. Supposons que les
fibres sont compactes. Quitte \`a conjuguer l'holonomie de la fibre
 au-dessus d'un point $x$ de $B$ par une translation $t_x$, on peut
supposer que les fibres ont m\^eme holonomie.
L'application $x\mapsto t_x$ peut \^etre choisie continue sur un
ouvert contenant $x$, on en d\'eduit de [G3] p. 178 que les
fibres au-dessus des points d'un voisinage de $p_B(x_0)$ sont
isomorphes entre elles, et par suite que toutes les fibres sont
isomorphes entre elles$\bullet$

\medskip

Le groupe $Aff(\hat
M,\hat\nabla_{\hat M})$ est muni  d'une topologie de Zariski d\'efinie
par les fonctions  $P\circ A_M$, o\`u $P$ est une
fonction polynomiale de $Aff({\R}^n)$.

\medskip

{\bf
Proposition 2.2.} {\it Soit $f$ un fibr\'e affine d'espace total
la vari\'et\'e affine  $(M,\nabla_M)$, si $\pi_1(F)$ est normal dans
$Z(\pi_1(M))$, l'adh\'erence de Zariski de $\pi_1(M)$ dans
$Aff(\hat M,\hat \nabla_{\hat M})$ alors pour tous points $\hat x$, et
$\hat y$  de $\hat M$ appartenant \`a la m\^eme orbite de
$Z(\pi_1(M))$, les fibres de $f(p_M(\hat x))$ et $f(p_M(\hat y))$
sont affinement isomorphes.}

\medskip

{\bf Preuve.} Tout
\'el\'ement $g$ de $Z(\pi_1(M))$ laisse stable $\hat{\cal F}^f$. On en
d\'eduit que $g$ induit un isomorphisme entre $\hat{\cal F}^f_{\hat
x}$ et $\hat{\cal F}^f_{g(\hat x)}$, et par suite un isomorphisme entre
les fibres respectives au-dessus de $f(p_M(\hat x))$ et $f(p_M(g(\hat
x)))$ car $g$ normalise $\pi_1(F)\bullet$

\medskip

{\bf
Remarque.}

 La proposition pr\'ec\'edente est vraie si $\pi_1(F)$ est
central dans $\pi_1(M)$. En effet, l'application  $Aff(\hat M,\hat
\nabla_{\hat M})\rightarrow Aff(\hat M,\hat\nabla_{\hat M})$ $\gamma
\mapsto \gamma\circ\alpha\circ\gamma^{-1}$ o\`u $\alpha$ est un
\'el\'ement de $\pi_1(F)$ est alg\'ebrique, comme elle est constante
sur $\pi_1(M)$, elle l'est aussi sur son adh\'erence de Zariski.

\medskip {\bf Corollaire 2.3.} {\it Si $(M,\nabla_M)$ est
compacte et compl\`ete, et $\pi_1(F)$ normal dans $Z(\pi_1(M))$, alors
toutes les fibres de $f$ sont affinement
isomorphes entre elles.}

\medskip

{\bf Preuve.} On sait que $Z(\pi_1(M))$
agit transitivement ([G-H2] th\'eor\`eme
2.6)$\bullet$

\medskip

Du corollaire pr\'ec\'edent on d\'eduit
que si l'espace total d'un fibr\'e affine est le tore de dimension
$n$, $T^n$ muni d'une structure affine compl\`ete,  toutes ses fibres
sont isomorphes. Ceci n'est pas tout \`a fait gratuit car sur $T^2$, il
existe des structures affines compl\`etes, qui ont la m\^eme holonomie
lin\'eaire sans \^etre isomorphes (voir paragraphe
1).

\medskip

{\bf Th\'eor\`eme 2.4} {\it Supposons que
l'adh\'erence de Zariski de l'holonomie de l'espace total
$(M,\nabla_M)$ d'un fibr\'e affine agisse transitivement sur
${\R}^n$, $\pi_2(B)=1$, et ses fibres  sont des vari\'et\'es
affines compactes et compl\`etes ayant un groupe fondamental
nilpotent, alors toutes les fibres sont isomorphes en tant que
vari\'et\'es affines.}

\medskip

{\bf Preuve.} Soient $x$ un point
de ${\R}^m$,  $\gamma$ un \'el\'ement de $\pi_1(B)$, et
$\hat\gamma$ un \'el\'ement de $\pi_1(M)$ au-dessus de $\gamma.$
 $\hat\gamma$ induit par conjugaison un automorphisme de $\pi_1(F)$,
$\Phi_{\hat\gamma}$.  Pour tout point $y$ de ${\R}^m$, la
restriction de $\pi_1(F)$ \`a $y\times {\R}^l$, est conjugu\'ee \`a
la restriction de $\pi_1(F)$ \`a $x\times {\R}^l$ par un
automorphisme polynomial $P_{\hat\gamma xy}$ (voir [F-G-H] th\'eor\`eme
7.1, proposition 8.3) induisant $\Phi_{\hat\gamma}$ sur $\pi_1(F)$.
Les coefficients de l'automorphisme $P_{\hat\gamma xy}$ sont des
fonctions polynomiales de $y$, car $\pi_1(F)$ est nilpotent (voir
[F-G-H]). $P_{\hat\gamma x\hat\gamma(x)}$ peut \^etre choisi affine,
\'egal \`a $\hat \gamma$ quitte \`a remplacer $P_{\hat\gamma xy}$
par $P_{\hat\gamma xy}P_{\hat\gamma x\hat\gamma(x)}^{-1}\hat\gamma$.
Ceci implique que les coefficients des monomes de
degr\'es sup\'erieurs ou \'egals \`a $2$,  de $P_{\hat\gamma xy}$,
s'annulent en $\hat\gamma(x)$. Ils s'annulent aussi en $P_{\hat\gamma
xg(x)}$, o\`u $g$ est un \'el\'ement de $Z(\hat\gamma)$ la cloture de
Zariski du groupe engendr\'e par $\hat\gamma$, car ces
coefficients sont des fonctions polynomiales de $y$. On en d\'eduit
que les fibres au-dessus des points de $p_B(Z(\hat\gamma)x)$ sont
isomorphes \`a la fibre au-dessus de $x$.  Le groupe $\pi_1(B)$
est engendr\'e par un nombre fini
d'\'el\'ements, $\gamma_1,...,\gamma_p$. Le groupe  $Z$ engendr\'e par
$Z(\pi_1(F)\hat\gamma_1)$,
 $Z(\pi_1(F)\hat\gamma_2)$, ...,$Z(\pi_1(F)\hat\gamma_p)$
est Zariski ferm\'e voir [Bor] p.190), en r\'ep\'etant
l'argument pr\'ec\'edent, on  d\'eduit que les fibres de $p_B(Zx)$
sont isomorphes entre elles.   Il en r\'esulte que toutes les fibres
sont isomorphes entre elles, car $Z$ contient l'adh\'erence de Zariski
de l'holonomie de l'espace total qui agit transitivement sur ${\R}^n\bullet$

\medskip

L'ensemble $App({\R}^m,{\R}^l)$ est
muni d'une structure de $\pi_1(F)$ module d\'efini
par $\gamma(D)=B_{\gamma}D$. Pour tout \'el\'ement $\gamma_1$
de $\pi_1(B)$,
consid\'erons l'application $k_{C_{\gamma_1}}:\pi_1(F)\rightarrow
App({\R}^m ,{\R}^l)$,
$\gamma\mapsto C_{\gamma}(h_B(\gamma_1))-C_{\gamma}$ est
un $1-$cocycle pour la structure de $\pi_1(F)$ module
pr\'ec\'edente.

\medskip

{\bf Proposition
2.5.} {\it Supposons que pour tout \'el\'ement $\gamma$ de $\pi_1(B)$,
la classe de cohomologie  $[k_{C_{\gamma}}]$, de $k_{C_{\gamma}}$, est
nulle et $Z(\pi_1(B))$ agit transitivement sur ${\R}^m$, alors
toutes les fibres sont isomorphes entre elles.} 

\medskip

{\bf
Preuve.} Si $[k_{C_{\gamma}}]$ est nulle, alors pour tout \'el\'ement
$\gamma$ de $\pi_1(B)$, $r(\gamma)(x)=r(x)$, l'application
$Aff({\R}^m )\rightarrow H^1(\pi_1(F),{\R}^l)$, $g\mapsto
r(g(x))-r(x)$ est alg\'ebrique et s'annulle sur $\pi_1(B)$, elle
s'annulle aussi sur sa cloture de Zariski. Comme $Z(\pi_1(B))$ agit
transitivement, on en d\'eduit que $r$ est une constante. Dans ce cas
le fibr\'e affine a des trivialisations locales affines.
(Voir paragraphe 4)$\bullet$

\medskip

 {\bf 3.
D\'eformation des structures affines \`a holonomie lin\'eaire
fix\'ee.}

\medskip

Consid\'erons $G$, le sous-groupe du groupe
des automorphismes de $\pi_1(F)$ v\'erifiant: pour tout \'el\'ement
$g\in G$, il existe une application lin\'eaire inversible $B_g$  de ${\R}^l$
telle que  $$L(h_F)(g(\gamma))=B_g\circ
L(h_F)(\gamma)\circ B_g^{-1}.\leqno (3.1)$$ Soit $B'_g$ une
autre application lin\'eaire v\'erifiant $3.1$. L'application
lin\'eaire de $B_g\circ {B'_g}^{-1}$ commute avec l'holonomie
lin\'eaire. $B_g$ est d\'etermin\'ee modulo les \'el\'ements du groupe
$Aut(L(h_F))$ des automorphismes lin\'eaires qui commutent avec
$L(h_F)$. L'ensemble $\{ B_g,g\in G\}$ sera appel\'e  groupe de gauge
des structures affines de $F$ dont l'holonomie lin\'eaire est
$L(h_F)$.

 L'application $:g\mapsto B_g$ est
une repr\'esentation de groupes modulo  les \'el\'ements de
$Aut(L(h_F))$.

 A tout $1-$cocycle $c$ de $L(h_F)$ et tout
\'el\'ement $B_g$ de $G$, on associe le $1-$cocycle $g^*c$ d\'efini par
 $g^*c(\gamma)=B_g(c(g^{-1}(\gamma))$. On d\'efinit ainsi une action
de $G$ sur $H^1(\pi_1(F),{\R}^l)$.

L'ensemble $E(F,L(h_F))$, des
classes d'\'equivalences  des structures affines de $F$ \`a holonomie
lin\'eaire \'egale \`a $L(h_F)$, est le quotient d'un ouvert
de $H^1(\pi_1(F),{\R}^l)$ par l'action de $G$.

\medskip

{\bf
Proposition 3.1.} {\it S'il existe une structure alg\'ebrique
de $E(F,L(h_F))$ rendant alg\'ebrique la projection canonique $t$ de
$H^1(\pi_1(F),{\R}^l)$ dans $E(F,L(h_F))$,   et si l'adh\'erence
de Zariski de l'holonomie de $(B,\nabla_B)$ agit transitivement
sur ${\R}^m $, alors toutes les fibres sont isomorphes en tant que
vari\'et\'es affines.}

\medskip

{\bf Preuve.} Soient $x$ un
point de ${\R}^m $, et $Z(\pi_1(B))$ l'adh\'erence de Zariski
de $h_B(\pi_1(B))$ dans $Aff({\R}^m )$. L'application
$Z(\pi_1(B))\rightarrow E(F,L(h_F))$,
$\gamma\mapsto t(r(\gamma(x))$ est alg\'ebrique, car $r$
est alg\'ebrique. Cette application est constante, car elle est
constante sur $\pi_1(B)$ qui est Zariski dense dans $Z(\pi_1(B))$. On
en d\'eduit que toutes les fibres sont isomorphes car $Z(\pi_1(B))$
agit transitivement sur ${\R}^m \bullet$

\medskip

{\bf
Remarque.}

 Les fibres d'un fibr\'e affine dont la base est compacte et
compl\`ete, et $E(F,L(h_F))$ est munie d'une structure alg\'ebrique
rendant alg\'ebrique la projection de $H^1(\pi_1(F),{\R}^l)$
sur  lui, sont isomorphes car l'adh\'erence de Zariski  de
$h_B(\pi_1(B))$ dans $Aff({\R}^m )$ agit transitivement sur ${\R}^m $.

\medskip

 {\bf Etude locale
de $E(F,L(h_F))$.}

\medskip

Pour d\'eterminer 'l'espace tangent'
en un point de $E(F,L(h_F))$, il faut \'etudier les d\'eformations de
ses \'el\'ements \`a holonomie lin\'eaire fix\'ee. Dans la suite de
cette partie on va supposer que $F$ est compacte et l'image de son
holonomie est discr\`ete. Ceci entraine  (voir [G3]) que l'ensemble des
repr\'esentations d'holonomies des structures affines de $F$ est
un ouvert dans l'ensemble des repr\'esentations de $\pi_1(F)$ dans
$Aff({\R}^l)$.

 Soient $h$ l'holonomie d'une structure affine
de $F$ dont la partie lin\'eaire est $L(h_F)$,   et $h_t$ une famille
\`a un  param\`etre de repr\'esentations de $\pi_1(F)$ dans $Aff({\R}^l)$,
 \`a parties lin\'eaires \'egales \`a $L(h_F)$ telles que
$h_0=h$. Il existe une application affine $u(x)$ de ${\R}^l$
pour tout \'el\'ement $x$ de $\pi_1(F)$, telle qu'on a:
 $$ h_t(x)
=exp(tu(x) + 0(t^2))h(x).$$

 Posons
$u=(B(x),b(x))$,
 $0(t^2)=(L(0(t^2))(x),a(0(t^2))(x))$,
et $h(x)=(A_x,a_x)$,
la partie lin\'eaire de   $h_t(x)$
est $exp(tB(x)+L(0(t^2))(x))A_x= A_x.$ En d\'erivant par rapport \`a
$t$, on obtient que
$exp(tB(x)+L(0(t^2)))(B(x)+L(0'(t^2))(x))A_x=0$,
 et par
suite $B(x)=0$ et $L(0(t^2))=0.$

 Ecrivons maintenant le fait que
$h_t$ est une repr\'esentation: on a
 $h_t(xy)=h_t(x)h_t(y)$,

 ceci est
\'equivalent
\`a
 $$ \matrix{(A_xA_y,A_x(a_y)+a_x+tb(xy)+ta(0(t^2))(xy))=
\cr (A_x,a_x+tb(x)+ta(0(t^2))(x))(A_y,a_y+tb(y)+ ta(0(t^2))(y))}$$

 et
par suite on d\'eduit que l'application $\pi_1(F)\rightarrow {\R}^l$,
 $x\mapsto b(x)$ est un $1-$cocycle pour $L(h_F)$.

 L'espace
tangent \`a un \'el\'ement de l'ensemble des repr\'esentations de
$\pi_1(F)$ dans $Aff({\R}^l)$, \`a partie lin\'eaire \'egale
\`a $L(h_F)$, s'identifie \`a l'ensemble $Z^1(\pi_1(F),L(h_F))$ des
$1-$cocycles de $\pi_1(F)$ pour la repr\'esentation
$L(h_F)$.

 D\'eterminons maintenant la structure de l'espace
tangent' d'un \'el\'ement de $E(F,L(h_F))$. Pour ce faire consid\'erons
une d\'eformation $h_t$ triviale de $h$ \`a partie lin\'eaire
$L(h_F)$. Il existe une courbe \`a un param\`etre $g_t$
de transformations du groupe de gauge, telle que
 $h_t(x)=g_{-t}h(x)g_t$.

 L'application $g_t$ est associ\'ee \`a l'identit\'e car on
a suppos\'e que l'image de l'holonomie est discr\`ete.
Posons
$g_t=exp(tu+0(t^2)), u=(C,c).$

 Comme la partie lin\'eaire de
$h_t(x)$ est celle de $h(x)$, $C$ commute avec $L(h_F)$. En
d\'erivant $exp(-t(u+0(t^2)))h(x)exp(tu+0(t^2))$ par rapport \`a $t$,
on obtient  $(0,-C(a_x)-c+A_x(c))$.

 "L'espace tangent" de
$E(\pi_1(F),L(h_F))$ en $h$ s'identifie au quotient $T_h$
de $Z^1(\pi_1(F),L(h_F))$ par l'ensemble des  $x\mapsto
(-C(a_x)-c+A_x(c))$, ou au quotient de $H^1(\pi_1(F),L(h_F))$
par $\{x\rightarrow [C(a_x)]\}$.

\medskip

{\bf Proposition 3.2.
Rigidit\'e} {\it Si $T_{h_{x_0}}$, l'espace tangent en la fibre
passant par $x_0$ est nul, alors toutes les fibres sont isomorphes
entre elles.}

\medskip

{\bf Preuve.} En effet l'ensemble des
\'el\'ements $x$ de $B$ tels que la fibre au-dessus de $x$, est
isomorphe \`a celle au-dessus de $x_0$, est un ferm\'e.
Si $T_{h_{x_0}}=0$, il est aussi un ouvert, d'o\`u le
r\'esultat$\bullet$

\medskip

Donnons maintenant des exemples de
$E(F,L(h_F))$.

 Si $F$ est une vari\'et\'e compacte munie
d'une m\'etrique plate, d'apr\`es les th\'eor\`emes de Bieberbach,
$E(F,L(h_F))$ est un point.

\medskip

{\bf Proposition 3.3.} {\it
Si $(F,\nabla_F)$ est une vari\'et\'e affine compacte, radiante, telle
que $\pi_1(F)$ est nilpotent, alors $E(F,L(h_F))$ est un
point.}

\medskip

{\bf Preuve.} En effet si $\pi_1(M)$ est
nilpotent d'apr\`es [F-G-H] toute structure affine  voisine de
celle d\'efinie par $h_F$ est radiante, il r\'esulte de [G3] qu'elle
lui est isomorphe, car quitte \`a conjuguer par une translation "assez
petite", on peut supposer qu'elles ont m\^emes holonomies$\bullet$

\medskip

Consid\'erons maintenant le tore $T^2$ muni
d'une structure affine compl\`ete distincte de la structure
riemannienne plate standard. Elle s'obtient en conjuguant le groupe
$\Gamma_{s,t}$ $s$, $t$ $\in {\R}$ engendr\'e par
$\gamma_{1,s}(x,y)=(x+s,y)$
 et $$\gamma_{2,t}(x,y)=(x,y+t)$$

par une application  polynomiale $\phi_r(x,y)=(x+ry^2,y), r\in
{\R}$, on obtient le groupe $\Gamma_{r,s,t}$ engendr\'e par
$\alpha_{1,s}(x,y)=(x+s,y)$
   et
$\alpha_{2,s,t}(x,y)=(x+2rty+rt^2,y+t)$.

Si on \'etudie les
d\'eformations de la structure plate en fixant l'holonomie lin\'eaire,
alors  $rt$ est une constante $a$; $\alpha_{2,s,t}$ s'\'ecrit
alors $\alpha_{2,s,t}(x,y)=(x+2ay+at,y+t)$.

 Le r\'eel $a$ \'etant donn\'e, pour tout
couple $(t,t')$ de r\'eels non nuls,  un automorphisme polynomial qui
conjugue $\Gamma_{r,s,t}$ en $\Gamma_{r',s',t'}$ est $f'= \phi_{a\over
t'}m\phi_{-a\over t}$, o\`u $m$ est l'automorphisme de ${\R}^2$
d\'efinit par $m(x,y)=({s'\over s}x,{t'\over t}y).$ On a
$$f'(x,y)=({s'\over s}x+a({st'-s't\over
st^2})y^2,{t'\over
t}y).$$

\medskip

{\bf Proposition 3.4.} {\it Les  vari\'et\'es
affines quotients de ${\R}^2$ respectivement par
$\Gamma_{{a\over t},s,t}$ et $\Gamma_{{a\over t'},s',t'}$
sont affinement isomorphes si et seulement
si $st'-s't=0$.}

\medskip

{\bf Preuve.} Supposons que
$st'-s't\neq 0$ et les structures affines d\'efinies
par $\Gamma_{{a\over t},s,t}$ et $\Gamma_{{a\over t'},s',t'}$ sont
isomorphes, alors il existe une application affine $B'$ qui conjugue
$\Gamma_{{a\over t'},s',t'}$ en $\Gamma_{{a\over t},s,t}$. La partie
lin\'eaire de $B'$ est triangulaire sup\'erieure car elle
pr\'eserve le groupe engendr\'ee la partie lin\'eaire de $\alpha_{2,s,t}$.
 Aussi
sans restreindre la g\'en\'eralit\'e on peut supposer que
 $B'$ conjugue
$\alpha_{1,s'}$ en $\alpha_{1,s}$. L'automorphisme $g=B'\circ f'$ est
un automorphisme polynomial qui normalise $\Gamma_{r,s,t}$, l'action
induite fixant $\alpha_{1,s}$, et
quitte \`a remplacer $g$  par son carr\'e, on peut supposer
 qu'il existe un \'el\'ement $p$ de
${\Z}$ tel
que $g(\alpha_{2,s,t})=\alpha_{2,s,t}+p\alpha_{1,s,t}$.
L'automorphisme de groupe r\'esultant de l'action de $B'$
sur $\Gamma_{{a\over t},s,t}$ peut s'exprimer sous la forme de la
matrice $C'=\pmatrix{1&p\cr 0&1}$. L'automorphisme $C$ de ${\R}^2$
dont la matrice est $\pmatrix{1&{ps\over t}\cr 0&1}$,
normalise $\Gamma_{{a\over t},s,t}$ et induit $C'$
sur $\Gamma_{{a\over t},s,t}$. $C^{-1}g$ commute avec
$\Gamma_{r,s,t}$. Il est montr\'e dans [T2] que l'ensemble
des automorphismes polynomiaux qui commute avec l'holonomie d'une
vari\'et\'e affine compacte et compl\`ete agit librement. Or on sait
que l'ensemble des automorphismes affines qui commute avec l'holonomie
d'une structure affine compl\`ete du tore agit transitivement car une
structure affine compl\`ete de $T^2$,  se rel\`eve en une structure
associative de ${\R}^2$, sous-jacente \`a sa structure de groupe
commutatif. Il en r\'esulte que  $C^{-1}g$ est affine. Il y a
contradiction $\bullet$

\medskip

Comme  $s$ et $t$ ne sont pas nuls, on
en d\'eduit que  l'ensemble des orbites des repr\'esentations affines
du tore $T^2$, dont la partie  lin\'eaire est l'holonomie
lin\'eaire d'une structure affine (fix\'ee) compl\`ete distincte de la
structure plate standard sous l'action du groupe de gauge est le
projectif r\'eel  priv\'e de deux points. On remarque que la projection
de $H^1(\pi_1(F),L(h_F))-\{0\}$ dans le projectif r\'eel est
alg\'ebrique. On en d\'eduit une proposition, qui donne un
crit\`ere, pour que  les fibres  d'un fibr\'e affine isomorphes
\`a des  tores munis de structures affines compl\`etes soient
isomorphes.

\medskip

{\bf Proposition 3.5.} {\it Soit $f$ un
fibr\'e affine dont les fibres sont des tores munies d'une structure
affine compl\`ete, alors si $Z(\pi_1(B))$, l'adh\'erence de Zariski de
$\pi_1(B)$ dans $Aff(\hat B,\nabla_{\hat B})$ agit transivement
sur $D_B(\hat B)$, alors toutes les fibres sont isomorphes entre
elles.}

\medskip

{\bf Preuve.} Si $L(h_F)$ est la
repr\'esentation triviale, le r\'esultat provient du fait que
$E(T^2,L(h_F))$ est un point. L'application
$r':Z(\pi_1(M))\rightarrow E(T^2,L(h_F))$ qui \`a $g$ associe la
classe d'\'equivalence de $r(gx_0)$ dans $E(T^2,L(h_F))$ est
alg\'ebrique on conclut grace \`a 3.1$\bullet$

\medskip

L'ensemble
$Hom(\Gamma, L)$ des repr\'esentations d'un groupe $\Gamma$ dans un
groupe de Lie $L$ est \'etudi\'e par plusieurs auteurs,
notamment lorsque $\Gamma$ est le groupe fondamental d'une surface
compacte. Si $\Gamma$ est de pr\'esentation fini, $Hom(\Gamma,L)$ a
une structure de vari\'et\'e alg\'ebrique sur laquelle $L$ agit en
composant une repr\'esentation par un automorphisme int\'erieur, le
quotient $Hom(\Gamma,L)/L$ de $Hom(\Gamma,L)$ par $L$ n'est pas en
g\'en\'eral une vari\'et\'e alg\'ebrique [G3].

 Dans l'\'etude de
l'espace de Teichmuller d'une surface $S$ de genre $g>1$, on est
amen\'e \`a consid\'erer les repr\'esentations de $\pi_1(S)$ dans
$PSL(2,{\R})$.

Supposons que $F$ soit compacte. L'ensemble
des repr\'esentations d'holonomies des structures affines de $F$
s'identifie \`a un ouvert de $Hom(\pi_1(F),Aff({\R}^l))$ [G3], et
l'ensemble des structures affines de $F$ au quotient de cet ouvert par
$Aff({\R}^l)$. Si $Hom(\pi_1(F),Aff({\R}^l))/Aff({\R}^l)$
est une munie d'une structure de vari\'et\'e alg\'ebrique rendant
alg\'ebrique la projection de $Hom(\pi_1(F),Aff({\R}^l))$ dans
$Hom(\pi_1(F),Aff({\R}^l))/Aff({\R}^l)$, alors toutes les
fibres d'un fibr\'e affine dont les fibres sont diff\'eomorphes
\`a $F$ et de base compacte et compl\`ete sont isomorphes entre
elles.

Soit $Ad(\pi_1(F)):  \pi_1(F)\rightarrow aff({\R}^l)$,
la compostion de l'holonomie et  de la repr\'esentation adjointe de
$Aff({\R}^l)$, $H^0(\pi_1(F),aff({\R}^l))$ est l'ensemble
des champs affines de $(F,\nabla_F)$. Les \'el\'ements de
$H^1(\pi_1(F),aff({\R}^l))$ dont le cup produit est nul sont
les \'el\'ements tangents \`a $Hom(\pi_1(F),Aff({\R}^l))/Aff({\R}^l)$
 voir [G2],  Il ne nous est pas connu si les groupes
de cohomologies sup\'erieurs $H^*(\pi_1(F),aff({\R}^l))$ ont
une interpr\'etation g\'eom\'etrique.

\bigskip

{\bf 4.
Fibr\'es affines affinement
localement triviaux.}

\medskip {\bf D\'efinitions.}

 On dira
qu'un fibr\'e affine   est un fibr\'e affine   trivial, s'il est le
produit affine  de sa base par une de ses fibres.

 On dira qu'un
fibr\'e affine est affinement  localement trivial (f.a.l.t), si et
seulement si son relev\'e  sur le rev\^etement universel de sa base
est un fibr\'e affine trivial.

\medskip

{\bf Remarque.}

Si
le fibr\'e affine $f$ est un f.a.l.t alors il existe un syst\`eme de
coordonn\'ees affines telles que pour tout \'el\'ement $\gamma$ de
$\pi_1(M)$ pr\'eservant une fibre,
 $h_M(\gamma)(x,y)=(x,B_\gamma(y)+d_\gamma)$; ceci \'equivaut \`a
dire que $r$ est constante. R\'eciproquement si le groupe des
transformations affines du rev\^etement universel de la base
agit transitivement sur celle-ci et $r=0$, le fibr\'e affine est un
f.a.l.t.

 Dans la suite de cette partie, on supposera que
la d\'eveloppante de l'espace total est injective et $0$ appartient
\`a son image.

\medskip

\medskip

{\bf Proposition
4.1.} {\it Soit $f$ un fibr\'e affine d'espace total une vari\'et\'e
affine  compl\`ete,  dont l'holonomie lin\'eaire des fibres est
celle d'une structure compl\`ete du tore $T^2$ distincte de la
structure riemannienne plate, alors, $f$ est un f.a.l.t.
}

\medskip

{\bf Preuve.} Ecrivons  ${\R}^n = {\R}^m \times {\R}^2$, l'holonomie de
la structure affine des fibres
de $f$ est engendr\'ee par deux \'el\'ements $\gamma_1$ et $\gamma_2$
de la forme: $$ \gamma_1(z,x,y)= (z, x+ry+
\alpha_1(z), y+\beta_1(z)). $$
 $$ \gamma_2(z,x,y)=(z,x+\alpha_2(z),y+\beta_2(z)). $$ O\`u,
$z$ et $ (x,y)$ d\'esignent respectivement des \'el\'ements de
 ${\R}^m $ et ${\R}^2$, $\alpha_1$, $\alpha_2$, $\beta_1$ et
$\beta_2$ des formes affines de ${\R}^m $ et $r$ un r\'eel non nul.
Comme $\gamma_1$ et $\gamma_2$ commutent, on
a: $$ \gamma_1\circ\gamma_2 =(z,x+\alpha_2(z)+
r(y+\beta_2(z))+\alpha_1(z), y+\beta_1(z)+\beta_2(z))$$ $=$ $$
\gamma_2\circ\gamma_1=(z,x+ry+\alpha_1(z)+\alpha_2(z),
y+\beta_1(z)+\beta_2(z)). $$ On
en d\'eduit que
 $\alpha_2(z)+r\beta_2(z)+\alpha_1(z)=\alpha_1(z)+ \alpha_2(z)$. Ceci
implique que $\beta_2(z)=0$.  Comme l'action de $\gamma_2$ sur ${\R}^n$
 est libre, on en d\'eduit que $\alpha_2$ est une constante
$a_2$. Aussi $\beta_1(z)$ ne peut s'annuler en un point $z_0$, sinon
 la fibre au-dessus de $f(z_0)$   ne serait pas compacte. On en
d\'eduit que $\beta_1$ est une constante $b_1$.

 Conjuguons
maintenant par
$$\phi(z,x,y)=(z,x,y-{1\over
r}\alpha_1(z))$$ $$\phi^{-1}(z,x,y)=(z,x,y+{1\over
r}\alpha_1(z))$$

 $$\phi^{-1}\circ\gamma_2\circ\phi =
\gamma_2$$

$$\phi^{-1}\circ\gamma_1\circ\phi(z,x,y)= (z,x+ry,y+b_1)\bullet$$

\medskip

 La
proposition pr\'ec\'edente n'est pas valable si on ne suppose pas que
la base est compl\`ete comme le montre le paragraphe
$1$.

\medskip

 Soient $(F,\nabla_F)$ et $(B,\nabla_B)$
deux vari\'et\'es affines compactes, on va classifier tous les
fibr\'es affines affinement  localement triviaux de
base $(B,\nabla_B)$ et de fibre $(F,\nabla_F)$
tels que la d\'eveloppante de l'espace total est injective.

 Soit $(M,\nabla_M)$
l'espace total d'un tel fibr\'e.  On sait
que $\pi_1(F)$ est un sous-groupe normal de $\pi_1(M)$, c'est le noyau
de la restriction de $h_M$ \`a ${\R}^m $. Consid\'erons $\gamma$
et $\gamma_1$ deux \'el\'ements quelconques de $\pi_1(F)$ et
de $\pi_1(M)$, on
a:

$$h_M(\gamma)(x,y)=(x,B_{\gamma}(y)+d_{\gamma})$$

 et
$$ h_M(\gamma_1)(x,y)=(A_{\gamma_1}(x)+
a_{\gamma_1},B_{\gamma_1}(y)+C_{\gamma_1}(x)+ d_{\gamma_1}), $$ o\`u
$A_{\gamma_1}$   d\'esigne  un automorphisme de ${\R}^m $, $B_{\gamma}$ et
 $B_{\gamma_1}$ des automorphismes de ${\R}^l$, $C_{\gamma_1}$ une
  application
  lin\'eaire de ${\R}^m $ dans
${\R}^l$, $a_{\gamma_1}$ un \'el\'ement de  ${\R}^m $,
$d_{\gamma}$ et $d_{\gamma_1}$ des \'el\'ements de  ${\R}^l$.

$$\matrix{h_M(\gamma_1)^{-1}(x,y)=
(A_{\gamma_1}^{-1}(x)- A_{\gamma_1}^{-1}(a_{\gamma_1}),
B_{\gamma_1}^{-1}(y)\cr
-B_{\gamma_1}^{-1} (d_{\gamma_1})-B_{\gamma_1}^{-1}
C_{\gamma_1}(A_{\gamma_1}^{-1}(x)- A_{\gamma_1}^{-1}(a_{\gamma_1})))}$$

\medskip

Ecrivons
que $h_M(\pi_1(F))$ est un sous-groupe normal de $h_M(\pi_1(M))$, on
a

$$\matrix{h_M(\gamma_1^{-1})\circ h_M(\gamma)\circ
h_M(\gamma_1)(x,y)=\cr
 (x,B_{\gamma_1}^{-1}B_{\gamma}B_{\gamma_1}(y)+
B_{\gamma_1}^{-1}B_{\gamma}C_{\gamma_1}(x)+ B_{\gamma_1}^{-1}B_{\gamma}
(d_{\gamma_1})\cr
 +B_{\gamma_1}^{-1}(d_{\gamma})-
B_{\gamma_1}^{-1}C_{\gamma_1}(x)- B_{\gamma_1}^{-1}(d_{\gamma_1}))}$$ Ceci
implique que
 $$B_{\gamma_1}^{-1}B_{\gamma}
  C_{\gamma_1}(x)-B_{\gamma_1}^{-1}C_{\gamma_1}(x)
 =0,$$ et
par suite, $C_{\gamma_1}(x)\in H^{0}(\pi_1(F),{\R}^l)$.

 L'espace
vectoriel  $Applin({\R}^m ,H^0(\pi_1(F),{\R}^l))$, des
applications lin\'eaires de ${\R}^m $ dans $H^0(\pi_1(F),{\R}^l)$,
 est muni d'une structure de $\pi_1(B)$ module \`a gauche
d\'efinie par $\gamma_1'(D)=B_{\gamma_1}\circ D$, et d'une structure
de $\pi_1(B)$ module  \`a droite d\'efinie par $\gamma_1'(D)= D\circ
A_{\gamma_1}$. O\`u $\gamma_1$ est un \'el\'ement de
$\pi_1(M)$ au-dessus de $\gamma_1'$.

\medskip

Notons $T_F$ la
composante connexe de l'ensemble  des transformations affines
de $(F,\nabla_F)$ qui se rel\`event sur ${\R}^l$ en des
translations.  L'application lin\'eaire de $H^0(\pi_1(F),{\R}^l)$
$t \mapsto  B_{\gamma_1}t$ induit sur   $Applin({\R}^m ,T_F)$
une structure de $\pi_1(B)$ module \`a gauche.  La structure de module
\`a droite de  $Applin({\R}^m ,H^0(\pi_1(F),{\R}^l))$, induit
sur $Applin({\R}^m ,T_F)$ une structure de $\pi_1(B)$ module \`a
droite. Soit ${\Z}\pi_1(B)$ l'alg\`ebre de groupe de $\pi_1(B)$.
L'espace vectoriel $Applin({\R}^m ,T_F)$ est donc muni d'une
structure de  ${\Z}\pi_1(B)$
bi-module.

\medskip

Consid\'erons l'application
 $C':\pi_1(B)\rightarrow  Applin({\R}^m , T_F$),
 $\gamma_1'\mapsto p_F(C_{\gamma_1})$, o\`u
 $p_F(C_{\gamma_1}(x))$ est l'\'el\'ement de $T_F$ induit par
$C_{\gamma_1}(x)$. V\'erifions que $C'$ est bien d\'efinie: Pour
tout \'el\'ement $\gamma(x,y) = (x,B_{\gamma}(y)+d_{\gamma})$
de $\pi_1(F)$, on a $$\gamma\circ
\gamma_1(x,y)= (A_{\gamma_1}(x)+a_{\gamma_1},
B_{\gamma}B_{\gamma_1}(y)+B_{\gamma}
C_{\gamma_1}(x)+B_{\gamma}(d_{\gamma_1})+d_{\gamma}),$$
comme $C_{\gamma_1}(x)\in
H^0(\pi_1(F),{\R}^l)$, on en d\'eduit que $C'$ est bien
d\'efinie.

Soit
$\gamma_2(x,y) =(A_{\gamma_2}(x)+a_{\gamma_2}, B_{\gamma_2}(y)+
C_{\gamma_2}(x)+d_{\gamma_2})$ un
autre \'el\'ement de $h_M(\pi_1(M))$, $$\matrix{\gamma_1\circ
\gamma_2(x,y) = (A_{\gamma_1}A_{\gamma_2}(x)+
A_{\gamma_1}(a_{\gamma_2})+a_{\gamma_1},\cr B_{\gamma_1}B_{\gamma_2}(y)+
B_{\gamma_1}C_{\gamma_2}(x)+ B_{\gamma_1}(d_{\gamma_2})+
C_{\gamma_1}(A_{\gamma_2}(x)+a_{\gamma_2})+ d_{\gamma_1})}.$$ Soient
$\gamma_1'$ et $\gamma_2'$, les images respectives de $\gamma_1$ et
$\gamma_2$ dans $\pi_1(B)$. On a $C'(\gamma_1'\circ\gamma_2')=
 p_F(B_{\gamma_1}C_{\gamma_2}(x)+C_{\gamma_1} A_{\gamma_2}(x))$,on
en d\'eduit que $C'$ se prolonge en $1-$cocycle $C$ de  ${\Z}\pi_1(B)$
 module pour la cohomologie de  Hochschild d\'efinit par les
actions \`a droite et \`a gauche pr\'ec\'edentes. Notons $[C]$ sa
classe de cohomologie.

L'application
de
$$C_2^f:\pi_1(B)\times\pi_1(B)\longrightarrow
T_F$$
$$ (\gamma_1,\gamma_2)\longmapsto
 p_F(C_{\gamma_1}(A_{\gamma_2})-C_{\gamma_1}))$$
  est un $2-$cocycle pour l'action \`a gauche
de $\pi_1(B)$.

\medskip

{\bf Proposition 4.2.}
{\it Supposons que la classe de cohomologie  de Hochschild $[C]$
associ\'ee \`a $f$, qu'on vient de d\'efinir s'annulle,
  alors $f$ est affinement  isomorphe \`a un fibr\'e affine
plat.}

\medskip

{\bf Preuve.} Supposons que $[C]=0$, alors il
existe une application $D'$ appartenant \`a $Applin({\R}^m ,T_F)$
telle que $C'(\gamma_1') = B_{\gamma_1}D'-D'A_{\gamma_1}$.  On en
d\'eduit que   $C_{\gamma_1} = B_{\gamma_1}\circ D- D\circ
A_{\gamma_1} +E$,  o\`u $D$ et $E$ sont des \'el\'ements de
 $Applin({\R}^m ,H^0(\pi_1(F),{\R}^l))$, tels que $p_F(D)=D'$
et $E(x)$ se projette en l'application identit\'e de $F$.
Comme $h_M(\pi_1(F))$ est discret car la d\'eveloppante de
$(M,\nabla_M)$ est injective, il en r\'esulte que $E=0$. En conjuguant
l'holonomie de $(M,\nabla_M)$ par $(x,y)\mapsto (x, y+D(x))$, on
en d\'eduit que $f$ est isomorphe \`a fibr\'e  affine
plat$\bullet$

\medskip

La donn\'ee de $[C]$ ne d\'etermine
pas compl\`etement la classe d'isomorphisme  du fibr\'e affine $f$,
car elle suppose la connaissance de l'action de $\pi_1(B)$ sur
$T_F$. Pour tout \'el\'ement $\gamma'_1$ de
$\pi_1(B)$, l'application
$(x,y)\mapsto (x,B_{\gamma_1}(y)+d_{\gamma_1})$   normalise
$h_M(\pi_1(F))$, elle d\'efinit donc un  \'el\'ement $g(\gamma'_1)$ de
$Aff(F,\nabla_F)$.

 On
a $$ (C'_{\gamma_1'\circ\gamma_2'})
g(\gamma_1'\circ\gamma_2')= (C'_{\gamma_1'}\gamma_2'-C'_{\gamma_1'})
 (C'_{\gamma_1'})g(\gamma_1')\circ (C'_{\gamma_2'})g(\gamma_2').$$
On dira que $g$ est compatible avec l'action  de
$\pi_1(B)$.

 Soit $Aff(F,\nabla_F)/T_F$ le quotient
de $Aff(F,\nabla_F)$ par $T_F$, $g(\gamma_1')$ se projette
sur $Aff(F,\nabla_F)/T_F$ en $\bar g(\gamma'_1)$. L'application
$\gamma_1'\mapsto  \bar g(\gamma'_1)$ est une
repr\'esentation.

\medskip

R\'eciproquement, si on se donne
une repr\'esentation $\bar g:
\pi_1(B)\rightarrow Aff(F,\nabla_F)/T_F$, pour tout \'el\'ement
 $\alpha$ de $\pi_1(B)$, consid\'erons un \'el\'ement $ g(\alpha)$ de
$Aff(F,\nabla_F)$ au-dessus de $\bar g(\alpha)$. L'application
$T_F\rightarrow T_F$, $t\mapsto g(\alpha)\circ
t\circ g(\alpha)^{-1}$ ne d\'epend que de  $\bar g(\alpha)$. On en
d\'eduit une repr\'esentation $\pi_1(B)\rightarrow
Aut(T_F)$, $t\mapsto B_{\alpha}(t).$ L'espace
vectoriel $Applin({\R}^m , T_F)$ est muni d'une structure de
 ${\Z}\pi_1(B)$ bimodule, la structure de module \`a gauche est  d\'efinie
par $(\alpha.C)(x) =B_{\alpha}C(x)$, et celle de module \`a droite est
d\'efinie par  $(C.\alpha)(x)= C(L(h_B(\alpha))(x))$. Si de plus on se
donne un $1-$cocycle  de Hochschild $C$ pour cette structure de
bimodule, telle que  l'application de  $g:\pi_1(B) \rightarrow
Aff(F,\nabla_F)$ au-dessus de $\bar g$
v\'erifie $$C_{\alpha_1\circ\alpha_2}
g(\alpha_1\circ\alpha_2)= (C_{\alpha_1}\alpha_2-C_{\alpha_1})
C_{\alpha_1}g(\alpha_1)\circ C_{\alpha_2}g(\alpha_2),$$
on d\'efinit un fibr\'e affine, en  quotientant  ${\R}^m \times
(F,\nabla_F)$ par l'action de $\pi_1(B)$ d\'efinie par $\alpha(x,y)=
(\alpha.x, C_{\alpha}(x)g(\alpha).(y))$.

\medskip

 Soient
$(F,\nabla_F)$, $(B,\nabla_B)$ deux vari\'et\'es affines et $f_1$
$f_2$ deux f.a.l.t de base $(B,\nabla_B)$ et de fibre $(F,\nabla_F)$.
A $f_1$ et $f_2$ on associe des repr\'esentations $\bar g_1$ et $\bar
g_2$ de $\pi_1(B)\rightarrow Aff(F,\nabla_F)/T_F$, induisant des
structures de ${\Z}\pi_1(B)$ bimodules sur $Applin({\R}^m ,T_F)$, deux
 $1-$cocycles $C_1$, $C_2$ pour la cohomologie de
Hochschild  de ces structures de bimodules, et enfin  deux applications
$g_1,g_2:\pi_1(B)\rightarrow Aff(F,\nabla_F)$ au-dessus de $\bar g_1$
et $\bar g_2$ telles que
 ${C_i}_{\alpha\circ\beta} g_{i}(\alpha\circ\beta)= 
({C_i}_{\alpha}(\beta)-{C_i}_\alpha){C_i}_\alpha g_i(\alpha){C_i}_\beta
g_i(\beta)$ $i\in\{1,2\}$.

\medskip

{\bf Th\'eor\`eme 4.3.} {\it
Soient $f_1$ et $f_2$ deux fibr\'es affines isomorphes, alors il existe
un automorphisme affine de $(B,\nabla_B)$ qui se rel\`eve sur $\hat
B\times (F,\nabla_F)$ en un automorphisme qui pr\'eserve les
fibres de $\hat f_1= \hat B\times (F,\nabla_F)$, conjugue
la repr\'esentation   $\bar g_1$ en $ \bar g_2$
et transforme $[C_1]$ en $[C_2]$. R\'eciroquement, s'il existe un
automorphisme affine de $(B,\nabla_B)$ qui se rel\`eve en
un automorphisme de  $\hat f_1$ pr\'eservant ses fibres,  conjuguant
$\bar g_1$ en $\bar g_2$, envoyant $[C_1]$ en $[C_2]$, alors il
existe un $1-$cocycle $c\in H^1(\pi_1(B),H^0(\pi_1(F),{\R}^l))$ tel que 
le fibr\'e affine dont l'holonomie est d\'efinie par
$\pi_1(M): \gamma \mapsto c(\gamma)h_{M_1}(\gamma)$ est isomorphe
\`a $f_2$. }

\medskip

On aura besoin du lemme
suivant:

\medskip

{\bf Lemme 4.4.} {\it Soient $C:{\R}^m \rightarrow H^0(\pi_1(F),{\R}^l)$
 une application lin\'eaire,
et $k=(B,d)$ (resp $k'=(A,a)$) un  \'el\'ement de $N(\pi_1(F))$, (resp
$N(\pi_1(B))$ le normalisateur de $\pi_1(F)$ (resp $\pi_1(B)$) dans
$Aff(\hat F,\hat\nabla_{\hat F})$, (resp $Aff(\hat B,\hat\nabla_{\hat
B}$)) l'automorphisme  $k_z:(x,y)\mapsto (A(x)+a,k(y)+C(z))$
conjugue $h_{M_1}$ en un groupe $h'_{M_1}$, il induit une nouvelle
repr\'esentation de $\pi_1(B)$ dans $T_F$. La classe de
cohomologie $[C_1]$ est transform\'ee en  une classe $[C'_1]$  pour
cette repr\'esentation qui ne d\'epend pas de $z$.}

\medskip

{\bf
Preuve.} Soit $\gamma$ un \'el\'ement de $h_{M_1}(\pi_1(M))$,
posons $$ \gamma(x,y)=(A_{\gamma}(x)+a_{\gamma},
B_{\gamma}(y)+C_{\gamma}(x)+d_{\gamma}),$$ on
a
 $$ \matrix{k_z^{-1}\gamma k_z=(k'^{-1}(A_{\gamma}(x),a_{\gamma})k',\cr
 B^{-1}(B_{\gamma}B(y)+B_{\gamma}C(z) +B_{\gamma}(d)+C_{\gamma}
 (k'(x))+d_{\gamma} -d-C(z)))}$$

La nouvelle repr\'esentation de $\pi_1(B)$ dans $T_F$ est induite par
 $$ h':\gamma\rightarrow B^{-1}B_1B$$ Le cocycle d\'efini par
$C_\gamma$ est transform\'e en celui d\'efini
par $$ (C'_\gamma)_z(x)= B^{-1}C_{\gamma}(A(x))\bullet$$

\medskip

{\bf
Preuve du th\'eor\`eme 4.3.} Soient $f_1$ et $f_2$ deux fibr\'es
affines isomorphes d'espaces total $(M_1,\nabla_{M_1})$, et
$(M_2,\nabla_{M_2})$. Il existe un automorphisme affine $k$ de
l'espace total de $\hat f_1$ respectant les feuilles de $\hat f_1$ tel
que   la restriction de $A_M(k)$ \`a ${\R}^m $ normalise
$h_B(\pi_1(B))$, et conjugue  $h_{M_1}$ en
$h_{M_2}$. Posons $$ A_M(k)(x,y)=(A(x)+a,B(y)+C(x)+d).$$

Pour
tout \'el\'ement $z$ de ${\R}^m $ tel qu'il existe $y$ dans ${\R}^l$
 v\'erifiant $(z,y)\in D_M(\hat M)$, l'automorphisme $A_M(k_z)$
  qui $(x,y)\mapsto (x,B(y)+C(z)+d)$ normalise l'holonomie de
$(F,\nabla_F)$, $A_M(k_z)$  d\'efinit  un \'el\'ement $k_z$ de
$Aff(F,\nabla_F)$ qui se projette sur $Aff(F,\nabla_F)/T_F$ en
un \'el\'ement $g$ ne d\'ependant pas de $z$. $g$ conjugue $\bar g_1$
en $\bar g_2$ et transforme $[C_1]$ en $[C_2]$.

R\'eciproquement,
supposons que les repr\'esentations $\bar g_1$ et $\bar g_2$ soient
conjugu\'ees par un automorphisme $g$ de ${\R}^m \times
(F,\nabla_F)$ laissant stable l'ensemble des fibres de $\hat B\times
(F,\nabla_F)$ et dont la restriction \`a ${\R}^m $ se restreint
sur $\hat B$ en un automorphisme qui normalise $\pi_1(B)$ (On a
suppos\'e que $D_M$ est injective). Supposons de plus, que $g$
transforme $[C_1]$ en $[C_2]$.

Conjuguons $g_1$ par $g$, on
obtient une application $g_1':\pi_1(B)\rightarrow Aff(F,\nabla_F)$
telle que pour tout \'el\'ement $\gamma$ de $\pi_1(B)$, on
a $g_1'(\gamma)= E(x)g_2(\gamma)$. Ou $E:{\R}^m \rightarrow T_F$
est une application affine. Comme $[C_1]$ est transform\'ee en $[C_2]$
par $g$, quitte \`a conjuguer $g_1'$ par une application
 $(x,y)\mapsto (x,y+C(x))$, o\`u $C:{\R}^m \rightarrow T_F$ est
lin\'eaire, on peut supposer que $E$ est une constante $c$.
L'application $E:\pi_1(B)\rightarrow T_F$ est un $\pi_1(B)$
$1-$cocycle$\bullet$

\medskip

{\bf Remarque.}

Etant
donn\'e un fibr\'e affine
affinement localement trivial au-dessus de $(B,\nabla_B)$ de fibre type
$(F,\nabla_F)$ dont la partie lin\'eaire coincide avec celle de
$h_{M_2}$, sa classe d'isomorphisme est d\'etermin\'ee par l'orbite de
$c$ (voir th\'eor\`eme pr\'ec\'edent) sous l'action du sous-groupe des
automorphismes affines de l'espace total qui pr\'eserve les fibres de
ce fibr\'e.

 Pour classifier les fibr\'es affines
affinement localement triviaux de base $(B,\nabla_B)$ et de fibre
$(F,\nabla_F)$, Il faut classifier

 - les classes de conjugaisons des
repr\'esentations de $\pi_1(B)\rightarrow Aff(F,\nabla_F)/T_F$,
cet ensemble peut  avoir une structure compliqu\'ee (voir
[G3]).

-  Une repr\'esentation
de $\pi_1(B)\rightarrow Aff(F,\nabla_F)/T_F$ \'etant donn\'ee, classifier:
  les
orbites  des $1-$cocycles de Hochschild de $Z\pi_1(B)$ \`a valeurs dans
$Applin({\R}^m , T_F)$ sous l'action de
$N(\pi_1(B))\times Aff(F,\nabla_F)$.

- Une
repr\'esentation de $\pi_1(B)\rightarrow Aff(F,\nabla_F)/T_F$ et
un $1-$cocycle \'etant donn\'es, classifier les applications $g$ compatibles
 avec
l'action de $\pi_1(B)$.

\medskip

{\bf Proposition 4.5.} {\it
Consid\'erons un fibr\'e affine caract\'eris\'e par $g, B_{\alpha}$ et
$C$, les fibr\'es caract\'eris\'es par $g', B_{\alpha} $
($B_{\alpha}$ ne varie pas)    sont de la forme $g'(\gamma)=
c({\gamma})g(\gamma)$, o\`u $c(\gamma)\in
H^1({\Z}\pi_1(B),App({\R}^m ,T_F))$. L'action \`a droite  de
$\pi_1(B)$ sur $Applin({\R}^m ,T_F)$ est induite par la
composition \`a droite avec l'holonomie de $(B,\nabla_B)$ et non par
la composition avec son holonomie lin\'eaire.}

\medskip

{\bf
Preuve.} On a $g'(\gamma) = c(\gamma)g(\gamma)$, o\`u $c(\gamma)$ est
un \'el\'ement de $App({\R} ,T_F)$. En \'ecrivant le fait que
 $$\matrix{c(\gamma_1\gamma_2)C_{\gamma_1\circ\gamma_2}
 g(\gamma_1\circ \gamma_2)
= ((C_{\gamma_1}+c(\gamma_1))(h_{M_2}(\gamma_2))-\cr
(C_{\gamma_1}+c(\gamma_1)) c(\gamma_1)C_{\gamma_1}g(\gamma_1)\circ
c(\gamma_2)C_{\gamma_2}g(\gamma_2),}$$
on obtient que  $c(\gamma_1 \gamma_2) =
B_{\gamma_1}c(\gamma_2)+ c(\gamma_1)\gamma_2$ d'o\`u le
r\'esultat$\bullet$

\medskip

{\bf Remarque.}

 Les cocycles de Hochschild
utilis\'es au  r\'esultat pr\'ec\'edent sont des applications affines
contrairement \`a ceux utilis\'es plus t\^ot qui sont des
applications lin\'eaires.

\medskip

On va faire
correspondre \`a nos $1-$cocyles de Hochschild des $1-$cocycles de
groupes.

A la repr\'esentation
$\pi_1(B)\rightarrow Aut(H^0(\pi_1(F),{\R}^l))$, on peut
associer un fibr\'e plat $BH^0$ qui est le quotient de $\hat B\times
H^0(\pi_1(F),{\R}^l)$ par l'action de $\pi_1(B)$ d\'efinie
par $(x,y)\mapsto (\gamma(x),B_{\gamma}(y))$. Il existe sur
$\Omega(B,H^0(\pi_1(F),{\R}^l))$, l'ensemble des formes d\'efinie
sur $B$ et \`a valeurs dans $BH^0$, une d\'erivation $d_1$ associ\'ee
\`a la structure plate de $BH^0$. L'ensemble $PBH^0$, des $1-$formes
$d_1$ parall\`eles est un faisceau localement constant au dessus de
$B$. Il s'identifie donc \`a un fibr\'e vectoriel au-dessus de $B$, qui
est le quotient de $\hat B\times Applin({\R}^m ,H^0(\pi_1(F),{\R}^l))$ par
 l'action de $\pi_1(B)$ d\'efinie par $(x,D)\mapsto
(\gamma(x), B_{\gamma}DA_{\gamma}^{-1})$.

Il existe
sur $\Omega(B,PBH^0)$ une d\'erivation $d_2$ d\'efinie par la
structure plate de $PBH^0$. Notons $B':\pi_1(B)\rightarrow
Aut(Applin({\R}^m ,H^0(\pi_1(F),{\R}^l)))$ la repr\'esentation
d\'efinie par $\gamma(D)= B_{\gamma}DA_{\gamma}^{-1}$. Les
$1-$cocycles pour cette repr\'esentation sont les
applications $C':\pi_1(B)\rightarrow
Applin({\R}^m ,H^0(\pi_1(F),{\R}^l))$
v\'erifiant $C'_{\gamma_1\gamma_2}=
B_{\gamma_1}C'_{\gamma_2}A_{\gamma_1}^{-1}+ C'_{\gamma_1}$.
A tout  $B'$ $1-$cocycle de groupe $C'$, on peut associer
 le ${\Z}\pi_1(B)$ $1-$cocycle de Hochschild  C d\'efini
par $C_{\gamma}=C'_{\gamma}A_{\gamma}$. L'application
 $C'\mapsto C$ est un isomorphisme d'espace vectoriel qui 
 passe au quotient en
un isomorphisme en cohomologie.

Supposons que $\hat B$ est
contractible. Les th\'eor\`emes d'isomorphismes de de Rham [Br]
montrent que $H^1(\pi_1(B),Applin({\R}^m ,{\R}^l))$ est
isomorphe \`a $H^1(B,PBH^0, d_2)$. Donc
\`a $C'$ associ\'e \` a un fibr\'e affine de base $(B,\nabla_B)$
et de fibe type $(F,\nabla_F)$ correspond un \'el\'ement $\hat C$
de $H^1(B,PBH^0,d_2)$. Aussi l'ensemble des $1-$ formes
$d_2$ ferm\'ees d\'efinies sur $B$ et \`a valeurs dans $PBH^0$
est isomorphe  \`a la somme de l'ensemble $\Omega^2_s(B,BH^0)$
des $2-$formes sym\'etriques $d_2$ ferm\'ees (c'est \`a dire $d_2$
parall\`eles), avec celles de l'ensemble des $2-$formes altern\'ees
$d_1$ ferm\'ees \`a valeurs dans $H^0(\pi_1(F),{\R}^l)$.
Donc \`a $\hat C$ il correspond
une deux forme altern\'ee $d_1$ ferm\'ee $\hat D$, et
une forme sym\'etrique parall\`ele
$D'$. Le
groupe $H^2(B,BH^0)$,  est isomorphe
\`a $H^2(\pi_1(B),H^0(\pi_1(F),{\R}^l))$ car $\hat B$ est
contractible. Donc il correspond
\`a la classe de $\hat D$ un \'el\'ement $[D]$ de
$H^2(\pi_1(B),H^0(\pi_1(F),{\R}^l))$.

On le th\'eor\`eme suivant:

\medskip

{\bf
Th\'eor\`eme 4.6.}
{\it Soient $f_1$ et $f_2$ deux fibr\'es
affines affinement localement triviaux au-dessus de $(B,\nabla_B)$  de
fibre type $(F,\nabla_F)$.
Supposons que $\hat B$ soit contractible. A $f_1$ et $f_2$ on associe
les repr\'esentations $\bar g_1,\bar
g_2:\pi_1(B)\rightarrow Aff(F,\nabla_F)/T_F$ et deux
\'el\'ements $[D_1],[D_2]$ et $D'_1, D'_2$
appartenant respectivement \`a $H^2(\pi_1(B),H^0(\pi_1(F),{\R}^l))$ et
 \`a $\Omega^2_s(B,BH^0)$ telles
que $d_2(D'_1)=d_2(D'_2)=0$ comme ci-dessus.

 Alors si $f_1$ et $f_2$ sont isomorphes, il existe
un automorphisme affine de $(B,\nabla_B)$ se relevant en un automorphisme
affine  du fibr\'e affine trivial $(\hat B,\hat\nabla_B)\times
(F,\nabla_F)$ de $Aff(F,\nabla_F)$ qui  conjugue $\bar g_1$ en  $\bar
g_2$ et transforme $([D_1],D'_1)$
en $([D_2],D'_2)$.

R\'eciproquement si $\bar g_1$ et $\bar g_2$
sont conjugu\'ees par un automorphisme affine de $(\hat B,\hat
\nabla_B)\times (F,\nabla_F)$ (se projettant sur $(B,\nabla_B)$) qui
transforme $([D_1],D'_1)$ en $([D_2],D'_2))$, alors il existe
un \'el\'ement $c\in H^1(\pi_1(B),H^0(\pi_1(F),T_F))$ tel que
le fibr\'e affine d\'efini par $g'(\gamma)=c(\gamma) g_1(\gamma)$
est isomorphe \`a $f_2$.}

\medskip

{\bf Proposition
4.7.}
{\it Supposons que l'espace total d'un fibr\'e affine $f$ soit
une vari\'et\'e affine $(M,\nabla_M)$ compacte et compl\`ete
et $dimH^1(\pi_1(F),{\R}^l)=1$, alors $f$ est un f.a.l.t. Si on
suppose de plus que la base  $(B,\nabla_B)$ est le tore $T^n $  muni de
sa structure riemannienne plate, et $dimH^1(\pi_1(M),{\R}^n)=1$
alors $f$ est un fibr\'e affine plat.}

\medskip

{\bf
Preuve.} L'application $r$ est  constante, sinon elle s'annullerait en
un point $x$ et la fibre au-dessus de $p_M(x)$ serait radiante, ceci
n'est pas possible car une vari\'et\'e affine compacte et compl\`ete ne
peut \^etre radiante. On en d\'eduit que $f$ est un f.a.l.t.

Si la
base est le tore $T^n$ muni de sa structure riemannienne plate
standard,  pour tout $x$ de ${\R}^m $,  $\gamma\mapsto
C_\gamma(x)$  est un $1-$cocycle  pour l'action de $\pi_1(B)$.
Notons $[C_\gamma(x)]$ sa classe de cohomologie.  La suite exacte de
Hochschild-Serre prouve que le groupe  $H^1(\pi_1(B),H^0(\pi_1(F),{\R}^n))$ 
s'injecte dans $H^1(\pi_1(M),{\R}^n)$,
l'obstruction radiante de $(M,\nabla_M)$ ne peut coincider avec un
\'el\'ement de la droite engendr\'ee par $[C_\gamma(x)]$ (s'il est non
nul), car l'action affine d'une vari\'et\'e affine compacte
et compl\`ete est irr\'eductible. On en d\'eduit que $[C_\gamma(x)]$
est nul, Et par  suite $f$ est  isomorphe \`a un fibr\'e  affine
 plat$\bullet$

\medskip

{\bf Proposition 4.8.}
{\it Soient $f_1$ et $f_2$
deux fibr\'es affines d'espace total compacts
respectifs $(M_1,\nabla_{M_1})$ et $(M_2,\nabla_{M_2})$, de m\^eme
base $(B,\nabla_B)$. On suppose en outre que $(B,\nabla_B)$ ne peut
 pas \^etre muni d'un feuilletage affine \`a feuilles compactes non trivial
  et
qu'il existe une fibre $F_0$  de $f_2$ qui ne peut pas \^etre l'espace
total d'un fibr\'e affine de base $(B,\nabla_B)$, alors tout
isomorphisme $g$ entre les vari\'et\'es affines $(M_1,\nabla_{M_1})$
et $(M_2,\nabla_{M_2})$ est un isomorphisme entre les fibr\'es
affines $f_1$ et $f_2$.}

\medskip

{\bf Preuve.} Soit $g$ un
isomorphisme entre  $(M_2,\nabla_{M_2})$ et  $(M_1,\nabla_{M_1})$,
 $f_1(g(F_0))$ n'est pas $B$ car $F_0$ n'est pas l'espace total d'un
fibr\'e affine au-dessus de $(B,\nabla_B)$. Les images des fibres
de $f_2$ par $f_1\circ g$ d\'efinissent sur $(B,\nabla_B)$ un
feuilletage affine \`a feuilles compactes.  Les  feuilles de
ce feuilletage affine sont forc\'ement  des points. On en d\'eduit
que $g$ envoie les fibres de $f_2$ sur celles de
$f_1\bullet$

\medskip

{\bf Proposition 4.9.}
{\it Supposons que
l'espace total $(M,\nabla_M)$ d'un fibr\'e affine $f$, soit une
vari\'et\'e affine compacte et compl\`ete, alors $Aff(M,\nabla_M)_0$
la composante connexe de $Aff(M,\nabla_M)$ pr\'eserve
$f$.}

\medskip

{\bf Preuve.}
Pour tout \'el\'ement $g$ de
$Aff(M,\nabla_M)_0$, il existe  un \'el\'ement $\hat g$ de
$Aff({\R}^n)$ au-dessus de $g$, qui commute avec $\pi_1(M)$. Soit
$0$ l'origine de ${\R}^n$. Le sous-groupe $\pi_1(F)$
de
$\pi_1(M)$ qui laisse stable $\hat{\cal F}^f_{0}$, laisse stable
$\hat g(\hat{\cal F}^f_{0})= \hat {\cal F}^{fg^{-1}}_{\hat g(0)}$, car
$\hat g$ commute avec $\pi_1(M)$, on en d\'eduit qu'il laisse stable
$\hat{\cal F}^{fg^{-1}}_{0}$. On sait que l'action affine de $\pi_1(F)$
sur  $\hat {\cal F}^f_0$ et $\hat {\cal F}^{fg^{-1}}_0$ est
irr\'eductible (voir [F-G-H] th\'eor\`eme 2.2). On en d\'eduit que les
deux sous-espaces vectoriels $\hat {\cal F}^f_0$ et $ \hat {\cal
F}^{fg^{-1}}_0$ coincident$\bullet$

\medskip {\bf
Arhitm\'eticit\'e de certaines classes d'isomorphismes des vari\'et\'es
affines.}

\medskip

Consid\'erons la vari\'et\'e de Hopf de
dimension $3$, $H_3$: c'est le quotient de ${\R}^3-\{0\}$ par une
homoth\'etie de rapport strictement positif
distinct de $1$. Cette vari\'et\'e affine
n'a pas de feuilletage affine \`a  feuilles   compactes. On va
d\'eterminer une classe de fibr\'es affines plats deux \`a deux
non isomorphes, dont la base de chaque \'el\'ement est $H_3$ et  la
fibre le tore de dimension $2$ munie de sa structure riemannienne
plate. En vertu de la proposition $4.8$, la classe d'isomorphisme
d'un tel fibr\'e affine plat d\'etermine celle de son espace
total.

Le groupe des transformations affines de $T^2$, $Aff(T^2)$
est isomorphe au produit semi-direct de $T^2$ par $Gl(2,{\Z})$.
Soient $\gamma$ le g\'en\'erateur de $\pi_1(H_3)$ et $B'_{\gamma}$ un
\'el\'ement de $Gl(2,{\Z})$ qui d\'etermine une repr\'esentation
$B_\pi$ de $\pi_1(H^3)$ dans $Aff(T^2)$. Soit $c:\pi_1(H_3)\rightarrow
T^2$, $\gamma\mapsto c_\gamma$, un $1-$ cocycle pour l'action de
$B_\pi$. L'application $\pi_1(H_3)\rightarrow
Aff(T^2)$, $\gamma\mapsto c_\gamma B_\gamma$ est
une repr\'esentation qui d\'efinit un fibr\'e affine plat au-dessus de
$H_3$. Si les valeurs propres de $B'_\gamma$ sont diff\'erentes de
$1$, alors $[c]$ est nul et le fibr\'e est isomorphe \`a celui
d\'efinit par $B_\pi$. Si les valeurs propres de $B'_\gamma$ sont $1$,
consid\'erons une base de $(e_1,e_2)$ de ${\R}^2$ telle
que $B'_\gamma(e_1)=e_1$, le cocycle $c_1$ d\'efinit par
$c_1(\gamma)=e_2$ n'est pas trivial. Soit $c_2$ un autre cocycle, si le
fibr\'e qu'il d\'efinit est isomorphe \`a celui d\'efinit par $c_1$,
alors il existe un \'el\'ement $D$ de $Aff(T^2)$ tel que
$D[c_1]=[c_2]$, car tout automorphisme affine de $H_3$ se rel\`eve
sur l'espace total, le fibr\'e \'etant plat. Si on choisit
 $c_2(\gamma)=he_2$ ou $h$ est un r\'eel non alg\'ebrique, on obtient
des fibr\'es affines non isomorphes et par suite des vari\'et\'es
affines qui ne sont pas isomorphes. Etudions maintenant une situation
plus g\'en\'erale.

A toute application $C:{\R}^3\rightarrow
{\R}^2$, on peut associer un $1-$cocycle de Hochschild de 
${\Z}\pi_1(H^3)$. Ce cocycle d\'etermine un fibr\'e affine en faisant le
quotient de ${\R}^3-\{0\}\times{\R}^2$, par le
groupe engendr\'e par  $\pi_1(T^2)$ (l'action de $\pi_1(T^2)$ se
fait par des translations qui se projettent en l'identit\'e
 sur ${\R}^3-\{0\}$), et le groupe engendr\'e par
l'application $(x,y)\mapsto (\lambda x,B_{\gamma}(y)+C(x))$, o\`u
$x$ et $y$ d\'esignent respectivement des \'el\'ements de ${\R}^3-\{0\}$
 et ${\R}^2$.  Ce cocycle est un cobord si et seulement
s'il existe une application affine $D$ de $App({\R}^3,{\R}^2)$
telle que $C=B_{\gamma}D-D\lambda=(B_{\gamma}-\lambda)D$. Choisisons
$C$ surjective, et supposons que $\lambda$ est une valeur propre de
$B_{\gamma}$, alors le cocycle $C'$ n'est pas exact car $C$
est surjective.

Si on fixe, $\lambda$ et $B_{\gamma}$, la
classe d'isomorphisme du fibr\'e affine obtenu est d\'etermin\'ee par
l'image de $[C']$ par $Gl({\R}^3)\times Gl(2,{\Z})$.
Ceci implique que les fibr\'es affines d\'etermin\'es par les $C$ sont
deux \`a deux isomorphes.

\bigskip

{\bf 5. Le cas
g\'en\'eral.}

\medskip

Soit $f$ un fibr\'e affine d'espace
total $(M,\nabla_M)$ et de base $(B,\nabla_B)$. On suppose ici que
$(M,\nabla_M)$ est compacte et compl\`ete, on a vu que dans ce cas
$\pi_1(M)$ est une extension de $\pi_1(B)$ par $\pi_1(F)$.  On va
aborder la classification des fibr\'es affines dans le cas
g\'en\'eral en utilisant la classification des extensions de
groupes.

Rappelons la classification des extensions de groupes (voir
[Mc] 124-137).  Soient     $Aut(\pi_1(F))$ le groupe des automorphismes
de $\pi_1(F)$, $Int(\pi_1(F))$ le groupe des
automorphismes int\'erieurs de $\pi_1(F)$, et $Out(\pi_1(F))$ le
quotient de $Aut(\pi_1(F))$ par $Int(\pi_1(F))$. Pour tout \'el\'ement
$\gamma$ de
$\pi_1(M)$, l'application $i^f_{\gamma}:\pi_1(F)\rightarrow
\pi_1(F)$, $\alpha\mapsto \gamma\circ\alpha\circ\gamma^{-1}$ passe
au quotient en une repr\'esentation $\Phi: \pi_1(B)\rightarrow
Out(\pi_1(F)).$

Soient $u:\pi_1(B)\rightarrow \pi_1(M)$ une section
telle que $u(1)=1$,  $v:\pi_1(B)\times \pi_1(B)\rightarrow
\pi_1(M)$ d\'efinie par  $v(x,y) = u(x)u(y)u(xy)^{-1}$ et
$\phi(x)\in \Phi(x)$. En appliquant la relation d'associativit\'e
\`a $u(x)u(y)u(z)$, on obtient
que
 $$i_{u(x)}v(y,z)v(x,yz)=v(x,y)v(xy,z). \leqno{(5.1)}$$ De
plus
$$\phi(x)\phi(y)=i_{v(x,y)}\phi(xy). \leqno{(5.2)}$$

Soient $K$
le centre de $\pi_1(B)$, $K$ est muni d'une structure de $\pi_1(B)$
module. Fixons $\Phi$, et supposons qu'il existe une extension de
$\pi_1(B)$ par $\pi_1(F)$. La donn\'ee d'un $2-$cocycle $w$ de
$\pi_1(B)$ \`a valeurs dans $K$ permet de d\'efinir une extension
de $\pi_1(B)$ par $\pi_1(F)$ en  posant
$(x_1,y_1)(x_2,y_2)= (x_1(\phi_1(y_1)(x_2))w(y_1,y_2),y_1y_2)$.
O\`u $x_1,x_2$ sont des \'el\'ements de $\pi_1(F)$ et $y_1,y_2$ des
\'el\'ements de $\pi_1(B)$. Pour $\Phi$ fix\'ee, les classes
d'\'equivalences d'extensions de $\pi_1(B)$ par $\pi_1(F)$ sont
donn\'ees par $H^2(\pi_1(B),K)$ s'il existe au moins une
extension.
La donn\'ee de $\Phi$ seul n'implique pas l'existence
d'une extension de $\pi_1(B)$ par $\pi_1(F)$ associ\'ee \`a $\Phi$.
L'obstruction \`a l'existence d'une telle extension est donn\'ee par
une famille d'\'el\'ements de $H^3(\pi_1(B),K)$.

\medskip
 La
classification des fibr\'es affines dans le cas g\'en\'eral va se faire
en fonction des donn\'ees suivantes:

(i) La base
$(B,\nabla_B)$.

(ii) La structure diff\'erentiable  de la fibre $F$
et une repr\'esentation d'holonomie lin\'eaire d'une de ses structures
affines

\medskip

Etant donn\'ee une extension $\pi_1(M)$ de $\pi_1(B)$ par
$\pi_1(F)$, qui est le groupe fondamental de l'espace total d'un
fibr\'e affine de base $(B,\nabla_B)$ dont l'holonomie lin\'eaire
des fibres est $L(h_F)$, on a vu que pour
tout \'el\'ement $\gamma\in\pi_1(M)$,
$$ h_M(\gamma)(x,y)=(A_{\gamma}(x)+a_{\gamma},
B_{\gamma}(y)+C_{\gamma}(x)+d_{\gamma})$$

On
a $$C_{\gamma\gamma'}=B_{\gamma}C_{\gamma'}+
C_{\gamma}(A_{\gamma'}),$$

En munissant $Applin({\R}^m ,{\R}^l)$ d'une structure de
$\pi_1(M)$ module \`a gauche d\'efinie
par $$ \gamma(C)=B_{\gamma}C, $$  et d'une structure de $\pi_1(M)$
module \`a droite d\'efinie par  $$ C\gamma=CA_{\gamma},$$  on
remarque que $C$ est un $1-Z\pi_1(M)$ cocycle pour l'homologie de
Hochschild d\'efinie par ses actions.

L'application
$$D:\pi_1(M)\rightarrow {\R}^l$$ $$ \gamma \mapsto
d_{\gamma} $$ est un $1-$cocycle pour l'action \`a gauche
de $\pi_1(M)$.

Le fibr\'e affine $f$ est affinement isomorphe \`a
un fibr\'e affine, affinement localement trivial si la restriction
de $C$ \`a $\pi_1(F)$ a une classe de cohomologie
nulle.

R\'eciproquement si on se donne une repr\'esentation
$B_\pi$ de $\pi_1(M)$ dans $Gl({\R}^l)$ dont la restriction
\`a $\pi_1(F)$ coincide avec $L(h_F)$, un $1-$cocycle $C$ pour
l'homologie de Hochschild pr\'ec\'edente et un $1-\pi_1(M)$ cocycle $D$
pour l' action \`a gauche de $\pi_1(M)$  tels que pour
tout \'el\'ement $x$ de ${\R}^m $, la
repr\'esentation $$\pi_1(F)\rightarrow Aff({\R}^l)$$
 $$ \gamma\mapsto (L(h_F)(\gamma), C(x)+d_{\gamma})$$ d\'efinit
une structure affine compl\`ete de $F$, alors on peut d\'efinir de
mani\`ere \'evidente une action affine de $\pi_1(M)$ sur ${\R}^n$,
dont l'espace quotient est l'espace total d'un fibr\'e affine au-dessus
de $(B,\nabla_B)$.

Le quotient de ${\R}^n$ par
$\pi_1(F)$ d\'efinit un fibr\'e affine au-dessus de ${\R}^m $ qu'on
note $\hat f$.

\medskip

{\bf Th\'eor\`eme 5.1.}
{\it Soient
$f_i,i\in \{1,2\}$ deux fibr\'es affines de base $(B,\nabla_B)$, dont
les fibres sont diff\'eomorphes et ont m\^eme holonomie lin\'eaire, et
 le groupe fondamental
de leur espace total respectif est $\pi_1(M)$, alors $f_1$ est
isomorphe \`a $f_2$ si et seulement s'il existe un isomorphisme entre
$\hat f_1$ et $\hat f_2$ au-dessus d'un automorphisme de
$(B,\nabla_B)$, dont un relev\'e sur ${\R}^n$
conjugue l'holonomie lin\'eaire de $(M_1,\nabla_{M_1})$ en
$(M_2,\nabla_{M_2})$,  transforme $[C_1]$ en $[C_2]$, et $[D_1]$ en
$[D_2]$.}

\medskip

{\bf Preuve.} Supposons que $f_1$ est
isomorphe \`a $f_2$, alors il existe une application affine $k$ qui
conjugue $h_{M_1}$ en $h_{M_2}$ en se projetant sur ${\R}^m $ en un
\'el\'ement de $N(\pi_1(B))$. $k$ conjugue aussi $h_{M_1}(\pi_1(F))$
en $h_{M_2}(\pi_1(F))$ ce qui entraine que $\hat f_1$ est isomorphe \`a
$\hat f_2$. Bien entendu, $[C_1]$ et $[D_1]$ sont respectivement 
transform\'es en
$[C_2]$ et $[D_2]$.

R\'eciproquement s'il existe un isomorphisme
entre $\hat f_1$ et $\hat f_2$ au-dessus d'un \'el\'ement de
$N(\pi_1(B))$, dont un relev\'e sur ${\R}^n$ conjugue
l'holonomie lin\'eaire de $(M_1,\nabla_{M_1})$ en celle de
$(M_2,\nabla_{M_2})$ et transforme $[C_1]$ (resp.$[D_1]$) en $[C_2]$
(resp. $[D_2]$), alors ce relev\'e conjugue $h_{M_1}$ en $h_{M_2}$
modulo une transformation affine de la forme $(x,y)\mapsto
(x,y+E(x)+e)\bullet$

\medskip

Soient $(B,\nabla_B)$ une vari\'et\'e
affine compacte et compl\`ete, et $F$ une vari\'et\'e diff\'erentiable
compacte  munie d'une structure affine compl\`ete dont l'holonomie
lin\'eaire est $L(h_F)$.

Pour classifier les fibr\'es affines
 d'espace total complet au-dessus de $(B,\nabla_B)$ de fibr
diff\'eomorphe \`a $F$ dont l'holonomie lin\'eaire  est $L(h_F)$ Il
faut:

- Classifier les extensions de $\pi_1(B)$ par $\pi_1(F)$,

 -  Une telle extension de $\pi_1(M)$ \'etant donn\'ee, classifier les
repr\'esentations lin\'eaires de $\pi_1(M)$ dans $Gl({\R}^m ,{\R}^l)$
 muni d'un $Z\pi_1(M))$ $1-$cocycle de Hochschild
$C\in Aff({\R}^m ,{\R}^l)$, et d'un $\pi_1(M)$ $1-$cocycle de
groupe $D$ \`a valeurs dans ${\R}^l$ tels que la repr\'esentation
$$\pi_1(F)\longrightarrow Aff({\R}^l)$$ $$ \gamma\longmapsto
(L(h_F)(\gamma),C(x)+D(\gamma)) $$

d\'efinisse une structure
affine.

\medskip

Les invariants utilis\'es pour la   classification dans
la situation g\'en\'erale sont plus extrins\`eques \`a la
g\'eom\'etrie affines que ceux utilis\'es dans le cas des fibr\'es
affines affinement localement triviaux.

Etant donn\'ee une
extension caract\'eris\'ee par $\Phi$ d\'efinissant un fibr\'e affine,
quelles sont les autres extensions d\'etermin\'ees par $\Phi$ donnant
lieu \`a des fibr\'es affines ?

On va aborder la classification
des fibr\'es affines de base $(B,\nabla_B)$ et de relev\'e $\hat f$
donn\'e.

Consid\'erons $\pi_1(M)$ une extension   de
$\pi_1(B)$ par $\pi_1(F)$, et supposons qu'il existe un fibr\'e affine
$f$, tel que le groupe fondamental de son espace total soit
$\pi_1(M)$. L'application $h_M(\gamma)$ se projette sur ${\R}^n/\pi_1(F)$
 en un automorphisme affine  $\bar h_M(\gamma')$.
L'application $\pi_1(B)\rightarrow Aut(\hat f)$, $\gamma'\mapsto
\bar h_M(\gamma')$ est une repr\'esentation. Soit $f'$ un autre
fibr\'e affine de base $(B,\nabla_B)$ tel que $\hat f=\hat f'$.
Posons $h_{M'}(\gamma)=h_M(\gamma)\circ   g(\gamma)$. L'application $
g(\gamma)$ se projette sur $\hat f$ en un automorphisme se projettant
en l'identit\'e sur ${\R}^m $. Pour tout \'el\'ement $\gamma$ de
$\pi_1(B)$, consid\'erons $ i^f(\gamma)$ l'application de
$Aut(\pi_1(F)) \rightarrow Aut(\pi_1(F))$ $\alpha\mapsto \bar
h_M(\gamma)\alpha\bar h_M(\gamma)^{-1}$.

En \'ecrivant le fait que
$\bar h'_M$ est une repr\'esentation, on obtient
 $$ g(\gamma_1\gamma_2)=i^f_{\gamma_2^{-1}}
 (g(\gamma_1)) g(\gamma_2). \leqno
(5.3)$$

R\'eciproquement la donn\'ee d'une application $
g:\pi_1(B)\rightarrow  Aut(\hat f)$ v\'erifiant $(5.3)$ permet
de d\'efinir un fibr\'e affine au-dessus de $(B,\nabla_B)$ en faisant
le quotient de $\hat f$ par le groupe engendr\'e par: $\bar
h_{M'}(\gamma)=\bar h_M(\gamma)g(\gamma).$

Sp\'ecialisons maintenant
au cas o\`u $i^f$ est fix\'ee, $ g(\gamma)$ commute avec
$\pi_1(F)$ puisque la fibre est compacte, $g(\gamma)$ appartient \`a
la composante connexe de $Aut(\hat f)$. La restriction
de $g(\gamma)$ \`a une fibre de $\hat f$ appartient \`a la composante connexe 
de son groupe de transformations affines voir [T3]). On a

$$
g(\gamma_1)g(\gamma_2)=g(\gamma_1\gamma_2).\leqno {(5.4)}
$$

R\'eciproquement la donn\'ee d'une repr\'esentation 
$g:\pi_1(B)\rightarrow Aut(\hat f)$ telle que $g(\gamma)$ appartient 
\` a composante connexe  de $Aut(\hat f)$ d\'etermine un fibr\'e affine
$f'$ au-dessus de $(B,\nabla_B)$ tel que $i_f=i_{f'}$.

\bigskip

Exprimons maintenant de mani\`ere plus g\'eom\'etrique l'application $r$
lorsque l'espace total du fibr\'e affine est compact et complet.
Notons $R(x)=r(x)-r(0)$. Fixons la base $(B,\nabla_B)$ et l'holonomie de 
la structure affine de la fibre au-dessus de $p_B(0)$.

Soit $\gamma$ un \'el\'ement de $\pi_1(M)$. Ecrivons 
$h_M(\gamma)(x,y)=(A_\gamma(x)+a_\gamma,B_\gamma(y)+C_\gamma(x)+d_\gamma)$.
L'\'el\'ement $\gamma$ induit par conjugaison un automorphisme $g_\gamma$
de $\pi_1(F)$. L'application $B_\gamma$ est un \'el\'ement du groupe de
gauge de l'holonomie lin\'eaire d'une fibre associ\'e \`a $g_\gamma$.
Si $\gamma$ appartient \`a $\pi_1(F)$, $g_\gamma$ laisse invariant $r(x)$
pour tout \'el\'ement $x$ de ${\R}^m$.

La repr\'esentation de $\pi_1(M)$ dans le groupe des automorphismes de 
$Gl(H^1(\pi_1(F),{\R}^l)^{\pi_1(F)})$ des automorphismes de
$H^1(\pi_1(F),{\R}^l)$, $\gamma\mapsto g_\gamma$, passe au quotient en une 
repr\'esentation $\bar g:\pi_1(B)\rightarrow 
Gl(H^1(\pi_1(F),{\R}^l)^{\pi_1(F)})$.
Il en r\'esulte un fibr\'e plat $H_{BF}$ au-dessus de $B$, qui est le quotient
de ${\R}^m\times H^1(\pi_1(F),{\R}^l)$ par l'action de $\pi_1(B)$
d\'efinie par $\gamma(x,y)=(\gamma(x),\bar g_\gamma(y))$.

Le fibr\'e $H_{BF}$ est munie d'une connexion plate canonique $\nabla_{BF}$,
\`a laquelle est associ\'ee la d\'eriv\'ee covariante $d_{BF}$. L'application
$R$ peut \^etre vue comme une $1-$forme $d_{BF}$ parall\`ele.

R\'eciproquement, consid\'erons une repr\'esentation $\bar g$ de 
$\pi_1(B)$ dans $Gl(H^1(\pi_1(F),{\R}^l)^{\pi_1(F)})$, telle que pour tout $\gamma$
appartenant \`a $\pi_1(B)$, il existe $B_\gamma$ appartenant au groupe
de gauge de l'holonomie lin\'eaire d'une fibre, dont l'action coincide avec
$g_\gamma$. La repr\'esentation $\bar g$ d\'efinit au-dessus de $B$ un fibr\'e
plat $H_{BF}$.

Pour se donner un fibr\'e affine au-dessus de $(B,\nabla_B)$ dont la structure
diff\'erentiable des fibres est $F$, et leur holonomie lin\'eaire est
$L(h_F)$, donnons nous d'abord une $1-$forme $d_{BF}$ parall\`ele
$R"$ \` a valeurs dans $H_{BF}$. La donn\'ee de $R"$ est \'equivalente 
\`a celle d'un $0-$cobord de groupe $R'$ de $\pi_1(B)$ \`a valeurs dans
$Applin({\R}^m,H^1(\pi_1(F),{\R}^l)^{\pi_1(F)})$ pour la structure de 
$\pi_1(B)$ module d\'efinie par $\gamma(D)=\bar g_\gamma DL(h_B)(\gamma^{-1})$.

Cette repr\'esentation d\'etermine un fibr\'e plat $H'_{BF}$ de fibre
type $Applin({\R}^m,H^1(\pi_1(F),{\R}^l))$. A $H'_{BF}$ on peut
associer la d\'erivation $d'_{BF}$. Le fibr\'e $H'_{BF}$ 
est isomorphe au faisceau 
localement constant des $1-$formes 
$d_{BF}'$ parall\`eles. Consid\'erons un \'el\'ement $R$ de 
$Applin({\R}^m,{\R}^l)$ au-dessus de $R'$. Supposons en outre que 
le $1-$cocycle de groupe $\gamma\mapsto R(\gamma)(x)+r(\gamma)(0)$ d\'efini
une structure affine de $F$. Le quotient de ${\R}^n$ par l'action de
$\pi_1(F)$ d\'efini pour tout \'el\'ement $\gamma$ de $\pi_1(F)$ par
$(x,y)\mapsto (x,B_\gamma(y)+R(\gamma)(x)+r(\gamma)(0))$, 
est un fibr\'e affine au-dessus de ${\R}^m$ tel que pour tout \'el\'ement
$\gamma_1$ de $\pi_1(B)$ et $x$ de ${\R}^m$, la fibre au-dessus de 
$p_B(x)$ est isomorphe \`a celle au-dessus de $p_B(\gamma(x))$.

\bigskip

\centerline{\bf Bibliographie.}

\bigskip
 
[Bo] Borel, A. Linear algebraic groups. W.A. Benjamin, Inc, New York-Amsterdam
1969.

[Br] Bredon, G. E. Sheaf theory. McGraw-HillBook Co., 1967.

[Ca] Carriere, Y. Autour de la conjecture de L. Markus sur les 
vari\'et\'es affines. Invent. Math. 95 (1989) 615-628.

[Ch] Charlap, L. S. Bieberbach groups and flat manifolds
Universitext. Springer-Verlag, 1986.

[F] Fried, D. Closed similarity affine manifolds.
Comment. Math. Helv. 55 (1980) 576-582.

[F-G] Fried, D. Goldman, W. Three-dimensional affine crystallographic
groups. Advances in Math. 47 (1983), 1-49.

[F-G-H] Fried, D. Goldman, W. Hirsch, M.
Affine manifolds with nilpotent holonomy. Comment. Math. Helv. 56
(1981) 487-523.

[G1] Goldman, W. Two examples of affine manifolds. Pacific J. Math.
94 (1981) 327-330.

[G2] Goldman, W. The symplectic nature of fundamental groups
of surfaces. Advances in in Math. 54 (1984) 200-225.

[G3] Goldman, W. Geometric structure on manifolds and varieties of
representations. 169-198, Contemp. Math., 74

[G-H1] Goldman, W. Hirsch, M. The radiance obstruction and parallel
forms on affine manifolds. Trans. Amer. Math. Soc. 286 (1984), 629-949.

[G-H2] Goldman, W. Hirsch, M. Affine manifolds and orbits of algebraic
groups. Trans. Amer. Math. Soc. 295 (1986), 175-198.

[God] Godbillon, C. Feuilletages. Etudes g\'eom\'etriques. Progress in
Mathematics, 98.

[K] Koszul, J-L. Vari\'et\'es localement plates et convexit\'e.
Osaka J. Math. (1965), 285-290.

[K] Koszul, J-L. D\'eformation des connexions localement plates.
Ann. Inst. Fourier 18 (1968), 103-114.

[Mc] Maclane, S. Homology. Springer-Verlag, 1963.

[Ma] Margulis, G. Complete affine locally flat manifolds with
a free fundamental group. J. Soviet. Math. 134 (1987), 129-134.

[Mi] Milnor, J. W. On fundamental groups of complete affinely flat
manifolds, Advances in Math. 25 (1977) 178-187.

[S-T] Sullivan, D. Thurston, W. Manifolds with canonical
coordinate charts: some examples. Enseign. Math 29 (1983), 15-25.

[T1] Tsemo, A. Th\`ese, Universit\'e de Montpellier II. 1999.

[T2] Tsemo, A. Automorphismes polynomiaux des vari\'et\'es affines.
C.R. Acad. Sci. Paris S\'erie I Math 329 (1999) 997-1002.

[T3] Tsemo, A. D\'ecomposition des vari\'et\'es affines.
Bull. Sci. Math. 125 (2001) 71-83.

[T4] Tsemo, A. Dynamique des vari\'et\'es affines. J. London Math. Soc. 63
(2001) 469-487.

[T5] Tsemo, A. From affine bundles to gerbes Preprint ICTP.

\end{document}